\documentclass[11pt,a4paper]{article}

 \usepackage{geometry}
\newgeometry{ hmargin={20mm,20mm}}

\usepackage{algorithm}
\usepackage{amsthm}
\usepackage{amsmath}
\usepackage{amssymb}
\usepackage{bm}
\usepackage{booktabs}
\usepackage{caption,subcaption}

\usepackage[numbers]{natbib}

\usepackage{tikz}
\usetikzlibrary{arrows.meta,
                chains,
                positioning,
                shapes.geometric
                }

\DeclareMathOperator*{\argmin}{argmin}
\DeclareMathOperator*{\minimize}{minimize}

\DeclareMathOperator{\diag}{diag}
\DeclareMathOperator{\Diag}{Diag}

\newtheorem{prop}{Proposition}

\DeclareMathOperator{\bA}{\bf A}

\DeclareMathOperator{\bD}{\bf D}
\DeclareMathOperator{\bd}{\bf d}
\DeclareMathOperator{\bE}{\bf E}

\DeclareMathOperator{\bI}{\bf I}

\DeclareMathOperator{\bl}{\bf l}
\DeclareMathOperator{\bM}{\bf M}
\DeclareMathOperator{\bL}{\bf L}
\DeclareMathOperator{\bp}{\bf p}
\DeclareMathOperator{\bQ}{\bf Q}

\DeclareMathOperator{\blr}{\bf r}

\DeclareMathOperator{\bu}{\bf u}

\DeclareMathOperator{\bv}{\bf v}

\DeclareMathOperator{\bw}{\bf w}

\DeclareMathOperator{\bx}{\bf x}

\DeclareMathOperator{\by}{\bf y}
\DeclareMathOperator{\bz}{\bf z}

\DeclareMathOperator{\betta}{\bm \eta}

\DeclareMathOperator{\bmu}{\bm \mu}

\DeclareMathOperator{\btheta}{\bm \theta}

\providecommand{\keywords}[1]{\textbf{\textit{Keywords:}} #1}

\begin{document}

\title{SCQPTH: an efficient differentiable splitting method for convex quadratic programming}
\author{Andrew Butler}
\maketitle

\begin{abstract}
We present SCQPTH: a differentiable first-order \textit{splitting} method for \textit{convex quadratic programs}. The SCQPTH framework is  based on the alternating direction method of multipliers (ADMM) and the software implementation is motivated by the state-of-the art solver OSQP: an operating splitting solver for convex quadratic programs (QPs). The SCQPTH software is made available as an open-source python package and contains many similar features including efficient reuse of matrix factorizations, infeasibility detection, automatic scaling and parameter selection. The forward pass algorithm performs operator splitting in the dimension of the original problem space and is therefore suitable for large scale QPs with $100-1000$ decision variables and thousands of constraints. Backpropagation is performed by implicit differentiation of the ADMM fixed-point mapping. Experiments demonstrate that for large scale QPs, SCQPTH can provide a $1\times - 10\times$ improvement in computational efficiency in comparison to existing differentiable QP solvers.  

\end{abstract}


\section{Introduction}\label{sec:intro}

Differentiable optimization layers enable the embedding of mathematical optimization programs as a custom layer in a larger neural network system. Notably, the OptNet layer, proposed by \citet{Amos2017}, implements a primal-dual interior point algorithm for constrained quadratic programs (QPs) with efficient backward implicit differentiation of the KKT optimality conditions. Similarly, \citet{Agrawal2019, Agrawal2020}, provide a differentiable optimization layer for general convex cone programming. Optimization programs are solved in the forward-pass by applying the alternating direction method of multipliers (ADMM) \citep{Gabay1976, Glow1975} to the corresponding homogeneous self-dual equations whereas the backward-pass performs implicit differentiation of the residual map. The aforementioned differentiable QP solvers are efficient for small scale problems, but can become computationally burdensome when the number of variables and/or constraints is large $( \gg 100)$.  More recently,  \citet{ButlerADMM}, provide an efficient ADMM implementation for solving box-constrained QPs, which for large scale QPs can provide up to an order of magnitude improvement in computational efficiency. The box-constrained QP layer, however, does not support more general linear inequality constraints. 

In this paper, we address this shortcoming and present SCQPTH: a differentiable first-order \textit{splitting} method for \textit{convex quadratic programs}.  The forward-pass implementation is based on ADMM and is highly motivated by the state-of-the art OSQP solver: an open-source, robust and scalable QP solver that has received strong adoption by both industry and academia \citep{Stellato2020}. SCQPTH contains many similar features including efficient reuse of matrix factorizations, infeasibility detection, automatic scaling and parameter selection. The OSQP solver, however, requires solving a system of equations of dimension $ \mathbb{R}^{(n+m ) \times (n+m)}$, which for large scale \textit{dense} QPs (i.e. with $n \gg 100$ and $m \geq  n$) can be computationally impractical. In contrast the SCQPTH software implementation follows the work of \citet{Boyd2011, Ghad2015} and others, and performs operator splitting in the dimension of the original problem space: $\bx \in \mathbb{R}^{n}$. The equivalent system of equations is therefore of dimension  $\mathbb{R}^{n \times n}$, and thus remains computational tractable for large scale dense QPs.  

The remainder of the paper is outlined as follows. In Section \ref{sec:method_admm_forward} we present the SCQPTH forward algorithm and discuss the relevant implementation features. In Section  \ref{sec:method_admm_backward} we follow  \citet{ButlerADMM} and present an efficient fixed-point implicit differentiation algorithm. In Section \ref{sec:results} we provide experimental  results comparing the computational performance of SCQPTH with the aforementioned differentiable QP solvers. We demonstrate that for large scale QPs, SCQPTH can provide a $1\times - 10\times$ improvement in computational efficiency.


\section{SCQPTH forward algorithm}  \label{sec:method_admm_forward}
We consider convex quadratic programming of the form:
\begin{equation}\label{eq:qp}
\begin{split}
\minimize \quad &   \frac{1}{2} \bx^T \bQ  \bx + \bx^T \bp\\
\text{subject to} \quad &   \bl \leq \bA \bx  \leq \bu, \\
\end{split}
\end{equation}
with decision variable $\bx \in \mathbb{R}^{n}$. The objective function is defined by a positive semi-definite matrix  $\bQ \in \mathbb{R}^{n \times n}$ and vector $\bp \in \mathbb{R}^{n}$. Linear equality and inequality constraints are defined by the matrix $\bA \in  \mathbb{R}^{m \times n}$ and vectors $\bl  \in \mathbb{R}^{m}$ and $\bu \in \mathbb{R}^{m}$. Quadratic programming occurs in many applications to statistics \citep{Tibshirani1996, Tikhonov1963},  machine-learning \citep{Vapnik1996, Ganti2011, Ho2015}, signal-processing \citep{Kim2008}, and finance \citep{Markowitz1952, Choueifaty2008, Bai2016}.

We solve program $\eqref{eq:qp}$ by applying the alternating direction method of multipliers (ADMM) algorithm.
Following \citet{Boyd2011, Ghad2015}, and others, we begin by introducing the auxiliary variable $\bz \in \mathbb{R}^m$ and recast program $\eqref{eq:qp}$ as:
\begin{equation}\label{eq:qp2}
\begin{split}
\minimize \quad &   \frac{1}{2} \bx^T \bQ  \bx + \bx^T \bp +  \mathbb{I}_{\bl \leq \bz \leq \bu}(\bz)\\
\text{subject to} \quad &   \bA \bx - \bz = 0. \\
\end{split}
\end{equation}
Following  \citet{Stellato2020}, the optimality conditions of program $\eqref{eq:qp2}$ are given by:
\begin{equation}\label{eq:qp_kkt}
\begin{split}
\bQ  \bx  + \bp + \bA^T\by & = 0\\
\bA \bx - \bz & =0\\
\bl \leq \bz & \leq \bu\\
\by_+^T (\bz - \bu) & = 0\\
\by_-^T (\bz - \bl)&  = 0,
\end{split}
\end{equation}
with Lagrange dual variable $\by$ and where $\by_+ = \max(\by, 0)$ and $\by_- = \min(\by, 0)$. The primal and dual residual of program $\eqref{eq:qp2}$ are therefore given as:
\begin{equation}\label{eq:qp_resid}
\blr_{\text{prim}} = \bA \bx - \bz \quad \text{and} \quad \blr_{\text{dual}} = \bQ  \bx  + \bp + \bA^T\by.
\end{equation}
Program $\eqref{eq:qp2}$ is separable in decision variables $\bx$ and $\bz$.  Let $\bmu = \rho^{-1}\by$, then the scaled ADMM iterations are as follows:

\begin{subequations} \label{eq:admm_qp}
\begin{align}
\bx^{k+1} & = \argmin_{\bx } \frac{1}{2} \bx^T \bQ \bx + \bx^T \bp + \frac{\rho}{2} \lVert \bA \bx - \bz^k + {\bmu}^k \rVert_2^2 \label{eq:admm_qp_x_1} \\
\bz^{k+1} & = \argmin_{\{ \bl \leq \bz \leq \bu \} } \frac{\rho}{2} \lVert \bA \bx^{k+1} - \bz + {\bmu}^k \rVert_2^2  \label{eq:admm_qp_z_1} \\
{\bmu}^{k+1} & = {\bmu}^{k} +  \bx^{k+1} - \bz^{k+1} \label{eq:admm_qp_mu_1}
\end{align}
\end{subequations}
Observe that Program $\eqref{eq:admm_qp_x_1}$ and Program $\eqref{eq:admm_qp_z_1}$ can be solved analytically and expressed in closed-form:

\begin{subequations} \label{eq:admm_qp_iter}
\begin{align}
\bx^{k+1} & = \alpha [\bQ + \rho \bA^T\bA]^{-1} (- \bp + \rho \bA^T (\bz^k - \bmu^k)) + (1-\alpha)\bx^{k}\label{eq:admm_qp_x} \\
\bz^{k+1} & = \Pi(\bA \bx^{k+1} + \bmu^k) \label{eq:admm_qp_z} \\
{\bmu}^{k+1} & = {\bmu}^{k} +  \bA\bx^{k+1} - \bz^{k+1} \label{eq:admm_qp_mu}
\end{align}
\end{subequations}
with relaxation parameter $0 < \alpha < 2$, step-size parameter  $\rho > 0$ and $\Pi$ denotes the projection operator onto the set $\{\bz \in \mathbb{R}^m \mid \bl \leq \bz \leq \bu\}$. We note that the matrix $[\bQ + \rho \bA^T\bA]$ is always invertible if $\bQ$ is positive definite or, without loss of generality, if $\bQ$ is positive semi-definite and $\bA$ has full row-rank. Indeed it is always possible to augment $\bA$ with the identity matrix and infinite bounds whenever $\bA$ is not full row rank. The per-iteration cost of the proposed ADMM algorithm is largely dictated by solving the linear system  $\eqref{eq:admm_qp_x}$. We note that this linear system is in general smaller than the KKT linear system in most interior point methods \citep{Nest1994}, homogenous self-dual embeddings \citep{Boyd2016}, and the OSQP ADMM implementation \citep{Stellato2020} and therefore remains computationally tractable when $m$ is large and either $\bQ$ and/or $\bA$ are dense. Indeed, the linear system of equations $\eqref{eq:admm_qp_x_1}$ is nearly identical to the indirect system proposed by \citet{Stellato2020} for large scale QPs; but does not require multiple conjugate gradient iterations. Lastly,  when $\rho$ is static then the algorithm requires only a single factorization of the matrix $[\bQ + \rho \bA^T\bA]$. We now discuss the forward solve implementation features, which in many cases is equivalent to the OSQP implementation and we refer the reader to \citet{Stellato2020} for more detail.

\subsection{Termination criteria}
When program $\eqref{eq:qp2}$ is strongly convex and the feasible set is nonempty and bounded then the ADMM iterations $\eqref{eq:admm_qp_iter}$ produces a convergence sequence such that:
\begin{equation}\label{eq:qp_conv}
\lim_{k \rightarrow \infty} \blr_{\text{prim}}^k = 0 \quad \text{and} \quad \lim_{k \rightarrow \infty} \blr_{\text{dual}}^k = 0.
\end{equation}
We refer the reader to \citet{Boyd2011} and \citet{Stellato2020} for proof. We define stopping tolerances as $\epsilon_\text{prim} > 0$ and  $\epsilon_\text{dual}  > 0$ for the primal and dual residuals, respectively. A reasonable stopping criteria is then given as:
\begin{equation}\label{eq:admm_stop}
|| \blr_{\text{prim}}^k ||_\infty \leq \epsilon_\text{prim} \quad \text{and} \quad || \blr_{\text{prim}}^k ||_\infty  \leq  \epsilon_\text{dual}.
\end{equation}
In practice we follow \citet{Boyd2011}  and define absolute and relative tolerances $\epsilon_\text{abs} > 0$ and $\epsilon_\text{rel} > 0$ and define stopping tolerances as :
\begin{equation}\label{eq:qp_tol}
\begin{split}
\epsilon_\text{prim} & = \epsilon_\text{abs} + \epsilon_\text{rel} \max( ||\bA \bx ||_\infty , ||\bz||_\infty )\\
\epsilon_\text{dual} & = \epsilon_\text{abs} + \epsilon_\text{rel}  \max( ||\bQ \bx ||_\infty , ||\bA^T\by||_\infty,  ||\bp||_\infty).
\end{split}
\end{equation}

\subsection{Infeasibility detection}
Alternatively, if program $\eqref{eq:qp2}$ is either unbounded or empty, then the ADMM iterations $\eqref{eq:admm_qp_iter}$ will terminate. Following \citet{Stellato2020}, program $\eqref{eq:qp2}$ is determined to be primal infeasible if the following conditions hold:
\begin{equation}\label{eq:prim_infeas}
\begin{split}
 ||\bA^T \Delta \by^k||_\infty  \leq \epsilon_\text{pinf}  ||\Delta \by^k||_\infty, \quad &  \bu^T(\Delta \by^k)_+ + \bl^T(\Delta \by^k)_-  \leq \epsilon_\text{pinf}  ||\Delta \by^k||_\infty
 \end{split}
\end{equation}
where  $ \epsilon_\text{pinf} > 0$ denotes the primal infeasibility tolerance and  $\Delta \by^k = \by^k - \by^{k-1}$.
Conversely, program $\eqref{eq:qp2}$ is determined to be dual infeasible if the following conditions hold:
\begin{equation}\label{eq:dual_infeas}
\begin{split}
||\bQ \Delta \bx ||_\infty  \leq \epsilon_\text{dinf}  ||\Delta \bx^k||_\infty, \quad &  ||\bp \Delta \bx ||_\infty  \leq \epsilon_\text{dinf}  ||\Delta \bx^k||_\infty \\
(\bA \Delta \bx)_i \quad & \begin{cases}
                \in [-\epsilon_\text{dinf}  ||\Delta \bx^k||_\infty, \epsilon_\text{dinf}  ||\Delta \bx^k||_\infty    ] & \text{if }  \bl_i, \bu_i \in \mathbb{R} \\
                \geq -\epsilon_\text{dinf}  ||\Delta \bx^k||_\infty& \text{if } \bu_i = \infty \\
                \leq \epsilon_\text{dinf}  ||\Delta \bx^k||_\infty & \text{if }\bl_i = \infty\
                \end{cases}
\end{split}
\end{equation}
with dual infeasibility tolerance $ \epsilon_\text{dinf} > 0$. 

\subsection{Scaling and parameter selection}
The ADMM algorithm is highly sensitive to the selection of the parameter $\rho$ and to the scale of the problem variables: $\bQ$, $\bp$, $\bA$, $\bl$ and $\bu$  \citep{Boyd2011} . Following \citet{Stellato2020} we propose a simplified scaling procedure that seeks to normalize the infinity norms of the matrices $\bQ$ and $\bA^T\bA$ and automatically selects $\rho$ based on their relative re-scaled norms. Specifically we let $0 \leq \beta \leq 1$ and define ${\bD = (1-\beta)\Diag(d_1, d_2, ..., d_n) + \beta \bar{d} }$ where $d_i = (|| \bQ_i ||_\infty)^{-\frac{1}{2}}$ and $\bar{d} = \text{mean}(d_i)$. The parameter $\beta$ therefore shrinks the norms, $d_i$, to the global mean, $\bar{d}$, which in practice can be helpful when the scale of the matrix $\bQ$ is not uniform. Similarly we let $\bE = \Diag(e_1, e_2, ..., e_m)$ where $e_i = || (\bA\bD)_i ||_\infty$. Program $\eqref{eq:qp}$ is then recast in the equivalent scaled form:

\begin{equation}\label{eq:qp_scaled}
\begin{split}
\minimize \quad &   \frac{1}{2} \bar{\bx}^T \bar{\bQ}  \bar{\bx} + \bar{\bx}^T \bar{\bp}\\
\text{subject to} \quad &   \bar{\bl} \leq \bar{\bA} \bar{\bx}  \leq \bar{\bu}, \\
\end{split}
\end{equation}
with scaled problem variables $\bar{\bQ} = \bD\bQ\bD$, $\bar{\bp} = \bD\bp$, $\bar{\bA} = \bE\bA\bD$, $\bar{\bl} = \bE\bl$ and $\bar{\bu} = \bE\bu$. The scaled decision variables are given as: $\bar{\bx} = \bD^{-1}\bx$, $\bar{\bz} = \bE\bz$ and $\bar{\by} = \bE^{-1} \by$. It is important to note that termination criteria and infeasibility detection are always performed with respect to the unscaled decision and problem variables. 

By default, we apply a heuristic initial parameterization for the step-size parameter $\rho$ given by the relative Frobenius norms: 
\begin{equation} \label{eq:rho_0}
\rho = \sqrt{\frac{m} {n}} \frac{||\bar{\bQ}||}{||\bar{\bA}^T\bar{\bA}||},
\end{equation}
where the scalar $\sqrt{m/n}$ corrects for the fact that $\bQ \in \mathbb{R}^{n \times n}$ and $\bA \in \mathbb{R}^{m  \times n}$. Thereafter, we follow \citet{Stellato2020} and adopt an adaptive parameter selection that updates $\rho$ based on the relative primal and dual errors, specifically: 
\begin{equation} \label{eq:rho_adapt}
\rho^{k+1} = \rho^k \sqrt{ \frac{ ||\blr_{\text{prim}}||_\infty / \max( ||\bA \bx ||_\infty , ||\bz||_\infty ) }{||\blr_{\text{dual}}||_\infty / \max( ||\bQ \bx ||_\infty , ||\bA^T\by||_\infty,  ||\bp||_\infty)} }. 
\end{equation}
In general it is computationally expensive to update $\rho$ at each iteration as this requires re-factorization of the matrix $[\bQ + \rho \bA^T\bA]$. Instead, we update $\rho$ when either $\rho^{k+1} / \rho^{k} > \tau$ or $\rho^{k} / \rho^{k+1} > \tau$, for some threshold parameter $\tau > 1$.


\section{SCQPTH backward algorithm}  \label{sec:method_admm_backward}
We now derive a method for efficiently computing the action of the Jacobian of the optimal solution, $\bx^*$, with respect to all of the QP input variables. We follow \citet{ButlerADMM} and recast the ADMM iterations in Equation $\eqref{eq:admm_qp_iter}$  as a fixed-point mapping of dimension $m$. We then apply the {\textbf{implicit function theorem}} \citep{Dontchev2009} in order to compute gradients of an arbitrary loss function $\ell  \colon \mathbb{R}^{n} \to \mathbb{R}$ with respect to the QP input variables. All proofs are available in the Appendix.

\begin{prop}\label{prop:admm_fixed_point}
Let $\bv^{k+1} = \bA \bx^{k+1} +\bmu^k$ and define $F \colon \mathbb{R}^{m} \times \mathbb{R}^{d_\theta} \to \mathbb{R}^{m}$. Then the ADMM iterations in Equation $\eqref{eq:admm_qp_iter}$  can be cast as a fixed-point iteration of the form $\bv = F(\bv,\btheta)$ given by:
\begin{equation}\label{eq:admm_fp}
\bv^{k+1} = \bA [\bQ + \rho \bA^T\bA]^{-1} (- \bp + \rho \bA^T (2\Pi(\bv^{k}) - \bv^{k})) + \bv^k - \Pi(\bv^{k})
\end{equation}
\end{prop}
The Jacobian, $\nabla_{\bv} F$, is therefore defined as:
\begin{equation}\label{eq:jacob_F}
\nabla_{\bv} F = \rho \bA [\bQ + \rho \bA^T\bA]^{-1} \bA^T (2D\Pi(\bv) - \bI) + \bI - D\Pi(\bv) 
\end{equation}
where $D\Pi(\bx)$ is the derivative of the projection operator, $\Pi$. The implicit function theorem therefore gives the desired Jacobian, $\nabla_{\btheta} \bv(\btheta)$, with respect to the the input variable $\btheta$:
\begin{equation}\label{eq:jacob_v}
\nabla_{\btheta} \bv(\btheta) = [\bI_{\bv} - \nabla_{\bv} F(\bv(\btheta),\btheta) ]^{-1} \nabla_{\btheta} F(\bv(\btheta),\btheta).
\end{equation}
From the definition of $\bv$ we have that the Jacobian $\nabla_{\btheta} \bx(\btheta) = \nabla_{\bv} \bx(\btheta) \nabla_{\btheta} \bv(\btheta)$ and where $$\nabla_{\bv} \bx(\btheta) = (\bA^T\bA)^{-1} \bA^TD\Pi(\bv).$$  In general forming the Jacobian $\nabla_{\btheta} \bx(\btheta)$ directly is inefficient, and instead we impute the  left  matrix-vector  product  of  the  Jacobian  with  the current gradient, $\frac{\partial \ell }{\partial \bx^*}$, as outlined below.

\begin{prop}\label{prop:admm_grads}
Let $\bd_{\bx}$ be defined as:
\begin{equation}\label{eq:grads_admm}
\bd_{\bx} = -[\bQ + \rho \bA^T\bA]^{-1} \bA^T [\bI_{\bv} - \nabla_{\bv} F(\bv(\btheta),\btheta) ]^{-T} [\nabla_{\bv} \bx(\btheta)]^T \frac{\partial \ell }{\partial \bx^*}
\end{equation}
Then the gradients of the loss function, $\ell$, with respect to input variables $\bQ$ and $\bp$ are given by:
\begin{equation}\label{eq:admm_partials}
\begin{aligned}
\frac{\partial \ell   }{\partial \bQ} & = \frac{1}{2} \Big(\bd_{\bx}  \bx^{*T} + \bx^* \bd_{\bx}^T \Big) & \qquad \frac{\partial \ell   }{\partial \bp} & = \bd_{\bx} \\
\end{aligned}
\end{equation}
\end{prop}
We approximate the gradients of the loss with respect to the constraint variables, $\bA$, $\bl$ and $\bu$ using the KKT optimality conditions.
\begin{prop}\label{prop:admm_grads_box}
Let $\bd_{\by}$ be defined as:
\begin{equation}
\bd_{\by} =  [\bA^T]^\dagger \Big( - \Big(\frac{\partial \ell }{\partial \bx^*} \Big)^T - \bQ  \bd_{\bx} \Big),
\end{equation}
where $ [\bA^T]^\dagger$ denotes the pseudo inverse of $\bA^T$. 
We define $\bd_{\by_-}$ and $\bd_{\by_+}$ as:
\begin{equation}\label{eq:dy_neg}
\bd_{\by_{-_j}}  = \begin{cases}
                \bd_{\by_j} / \by_{-_j} & \text{if }   \by_{-_j} \leq 0\\
                0 & \text{otherwise,} \\
                \end{cases} \qquad \text{and} \qquad
\bd_{\by_{+_j}} =  \begin{cases}
                \bd_{\by_j} / \by_{+_j} & \text{if }   \by_{+_j} \geq 0\\
                0 & \text{otherwise,} \\
                \end{cases}
\end{equation}
Then the gradients of the loss function, $\ell$, with respect to problem variables $\bA$, $\bl$ and $\bu$ are given by:
\begin{equation}\label{eq:admm_partials_box}
\frac{\partial \ell   }{\partial \bA} =  \bd_{\by} \bx^{*T} + \by^*  \bd_{\bx} ^T  \qquad  \frac{\partial \ell   }{\partial \bl}  = \diag(\by^*_-) \bd_{\by_-} \qquad  \frac{\partial \ell   }{\partial \bu} =\diag(\by^*_+) \bd_{\by_+}.
\end{equation}
\end{prop}

\section{Computational experiments}\label{sec:results}
We present three experiments that evaluates the computational efficiency of the SCQPTH implementation. Computational efficiency is measured by the median runtime required to execute the forward and backward algorithms. We compare against 3 alternative methods:
\begin{enumerate}
\item \textbf{LQP:} the box-constrained QP layer implementation proposed by \citet{ButlerADMM}. Implements ADMM  in the forward-pass and fixed-point implicit differentiation in the backward-pass.

\item \textbf{QPTH:} the OptNet layer implementation proposed by \citet{Amos2017}. Implements a primal-dual interior-point method in the forward-pass and efficient pre-factorized KKT implicit differentiation in the backward-pass.

\item \textbf{Cvxpylayers:} the general convex layer implementation proposed by \citet{Agrawal2019}. Implements the Splitting cone solver (SCS) \citep{Boyd2016} in the forward-pass and conic residual map implicit differentiation in the backward-pass.

\end{enumerate}

All experiments are conducted on an Apple Macbook Pro computer (2.6 GHz 6-Core Intel Core i7,32 GB 2667 MHz DDR3 RAM) running macOS `Monterey' and Python 3.9.

\subsection{Experiment 1: Random Box  Constrained QPs}\label{sec:results_1}
We generate random box constrained QPs of dimension: $n \in \{10, 25, 50, 100, 250, 500, 740, 1000\}$ and consider low $(1\mathrm{e}{-3})$ and high $(1\mathrm{e}{-5})$ absolute and relative stopping tolerances.  Problem variables are generated as follows. We set $\bQ =  \bL^T\bL + 0.01 \bI$ where $\bL \in \mathbb{R}^{n \times n}$ and entries $\bL_{ij} \sim \mathcal{N}(0, 1)$ with $50\%$ probability of being non-zero. Similarly we let $\bp_i \sim \mathcal{N}(0, 1)$,  $\bl_j \sim \mathcal{U}(-2,-1)$ and $\bu_j \sim \mathcal{U}(1,2)$. Experiment results are averaged over 10 independent trials with a batch size of 32.

\begin{figure}[h]
  \centering
  \begin{subfigure}[b]{0.22\linewidth}
    \includegraphics[width=\linewidth , trim={0mm 0cm 0cm 0cm},clip]{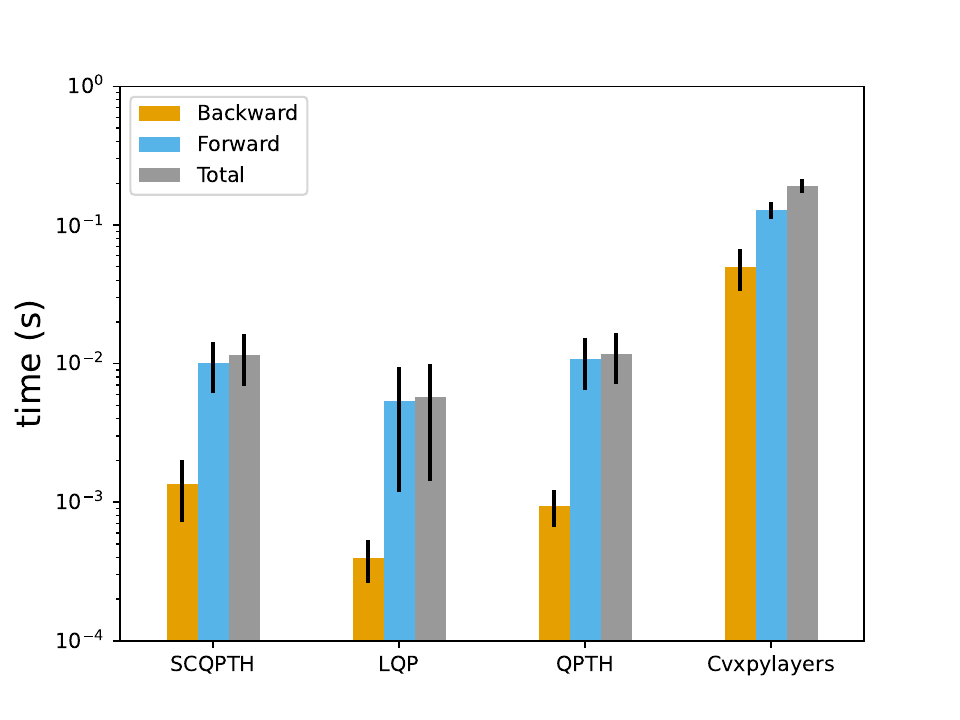}
    \caption{$n = 10$.}
  \end{subfigure}
  \begin{subfigure}[b]{0.22\linewidth}
    \includegraphics[width=\linewidth , trim={0mm 0cm 0cm 0cm},clip]{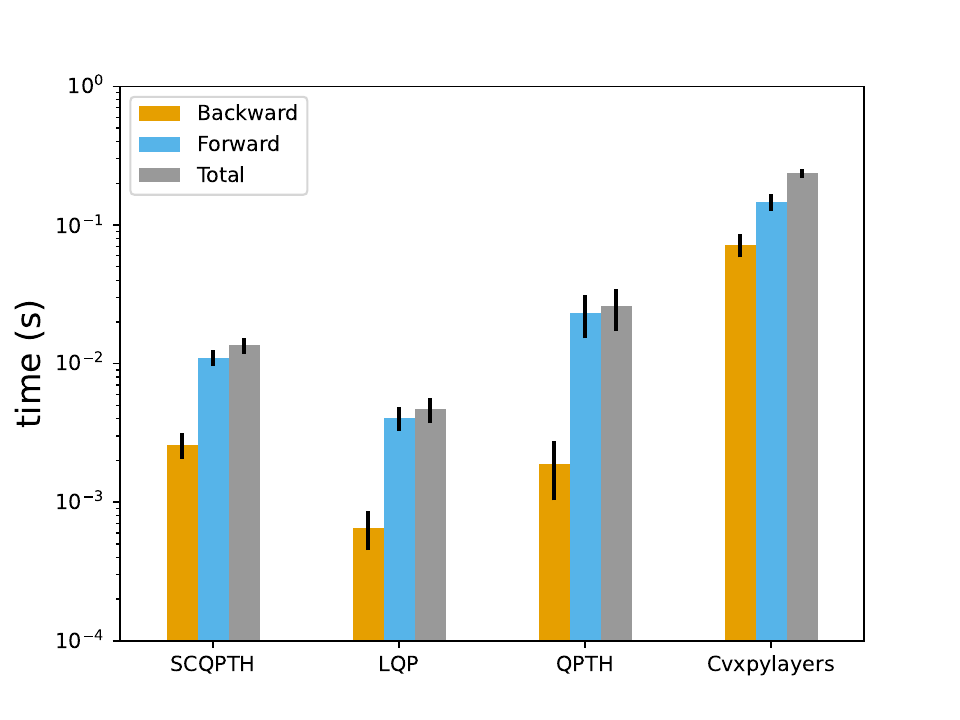}
    \caption{$n = 25$.}
  \end{subfigure}
    \begin{subfigure}[b]{0.22\linewidth}
   \includegraphics[width=\linewidth , trim={0mm 0cm 0cm 0cm},clip]{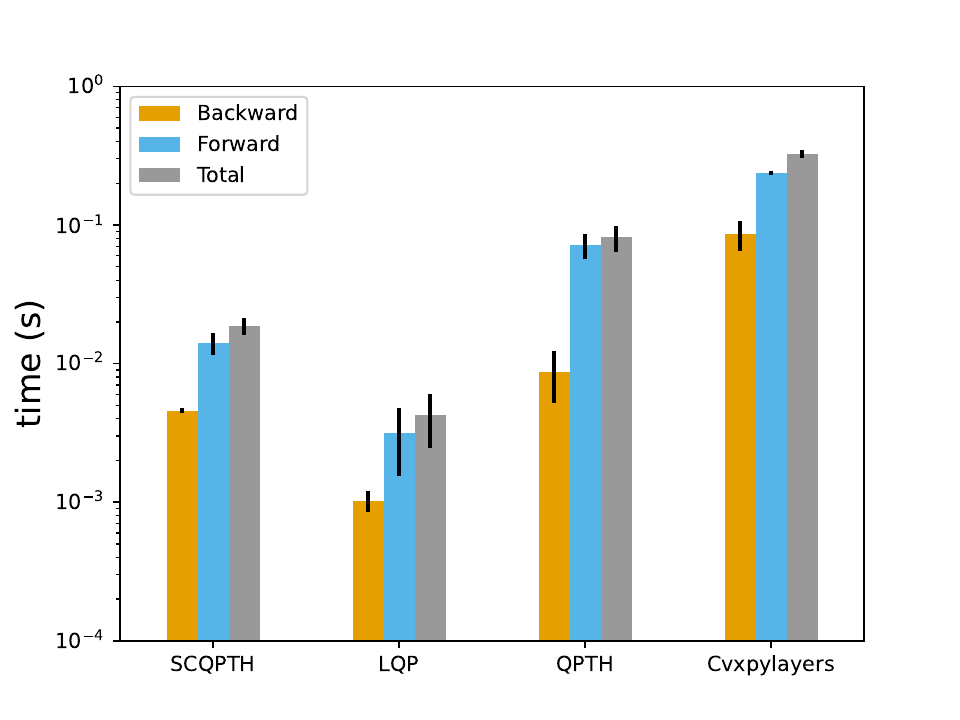}
    \caption{$n= 50$.}
  \end{subfigure}
  \begin{subfigure}[b]{0.22\linewidth}
    \includegraphics[width=\linewidth , trim={0mm 0cm 0cm 0cm},clip]{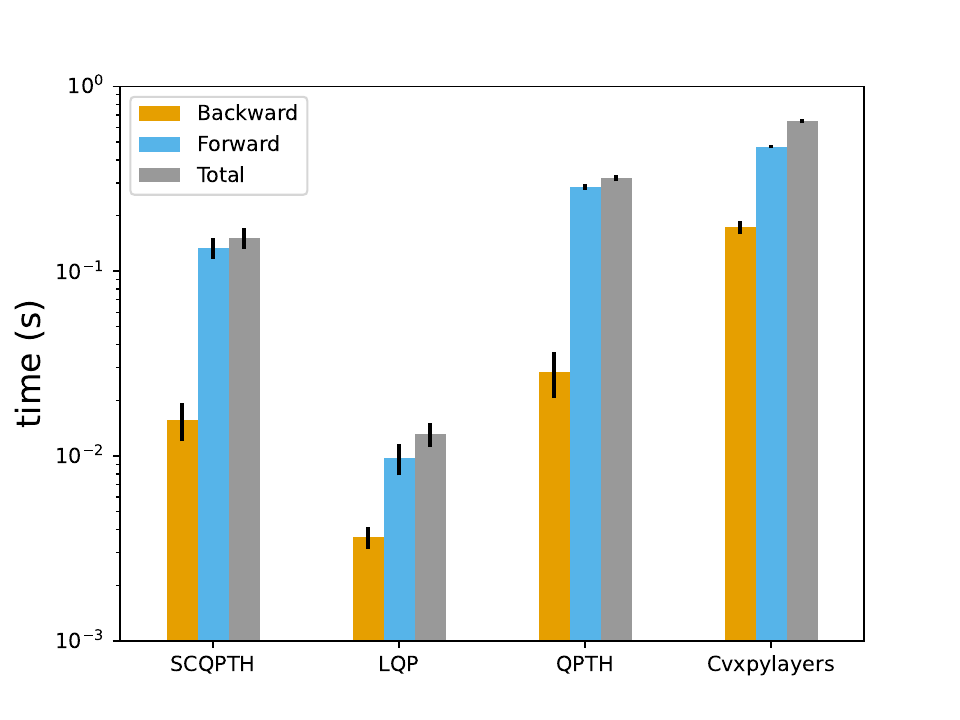}
    \caption{$n = 100$.}
  \end{subfigure}
  \begin{subfigure}[b]{0.22\linewidth}
    \includegraphics[width=\linewidth , trim={0mm 0cm 0cm 0cm},clip]{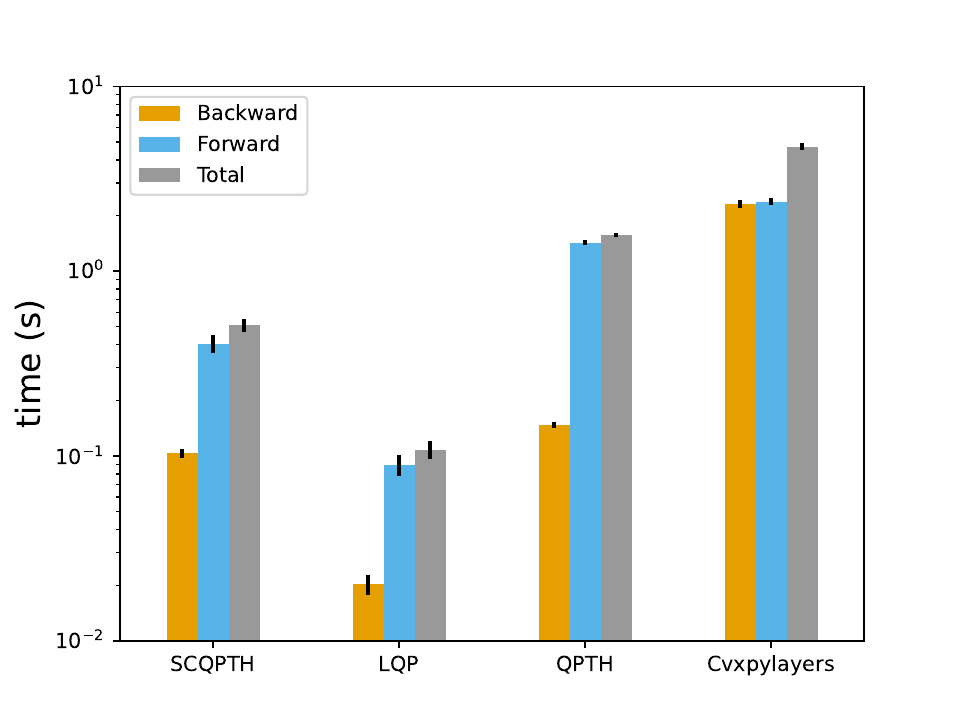}
    \caption{$n = 250$.}
  \end{subfigure}
  \begin{subfigure}[b]{0.22\linewidth}
    \includegraphics[width=\linewidth , trim={0mm 0cm 0cm 0cm},clip]{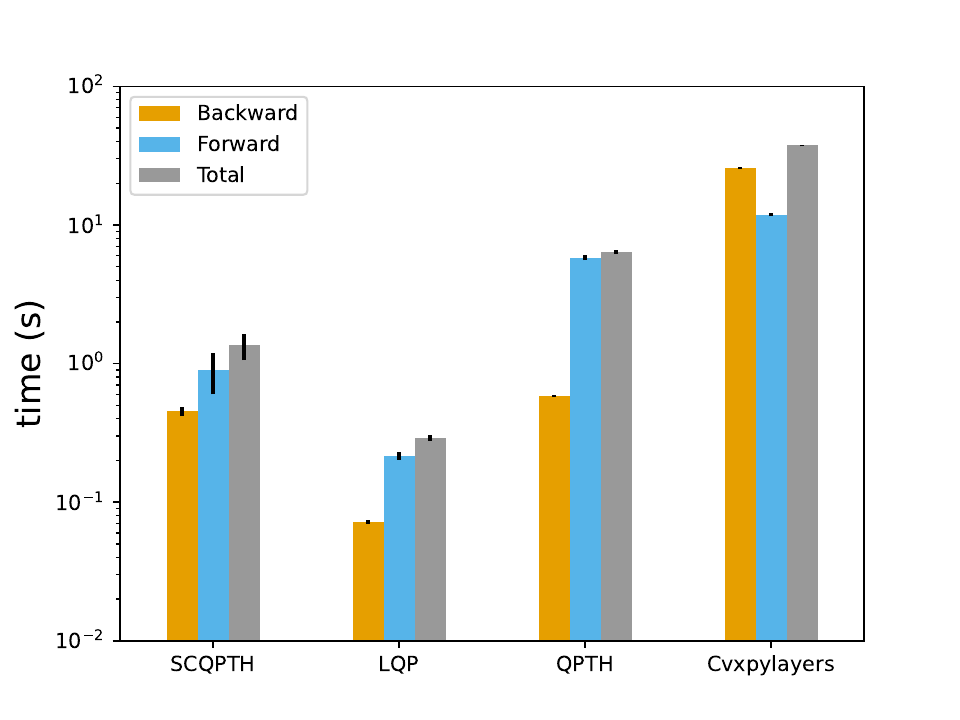}
    \caption{$n = 500$.}
  \end{subfigure}
  \begin{subfigure}[b]{0.22\linewidth}
    \includegraphics[width=\linewidth , trim={0mm 0cm 0cm 0cm},clip]{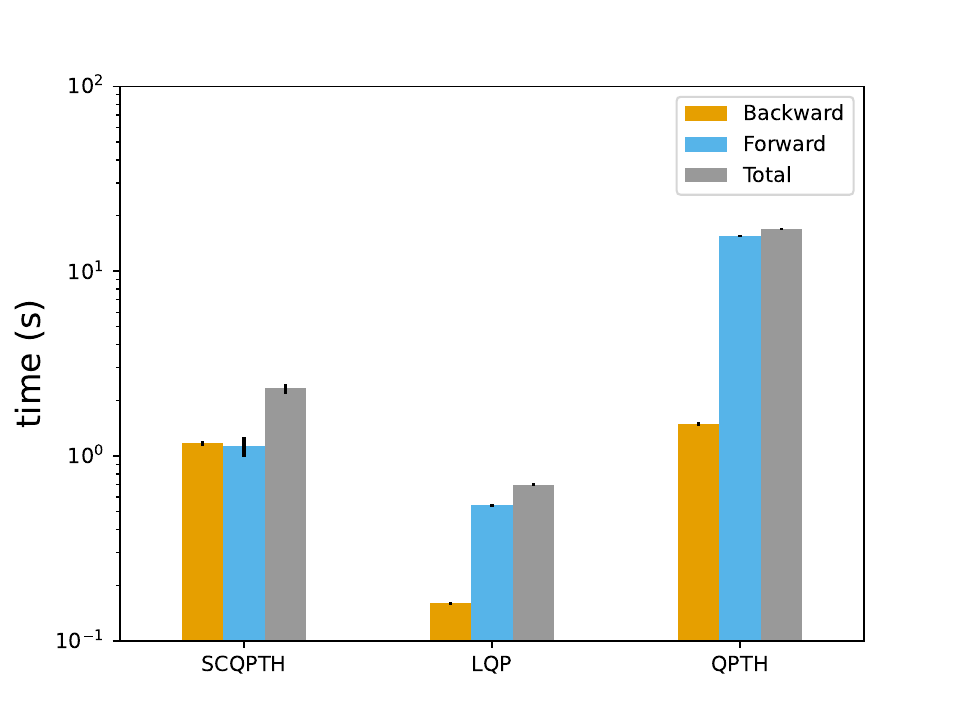}
    \caption{$n = 750$.}
  \end{subfigure}
  \begin{subfigure}[b]{0.22\linewidth}
    \includegraphics[width=\linewidth , trim={0mm 0cm 0cm 0cm},clip]{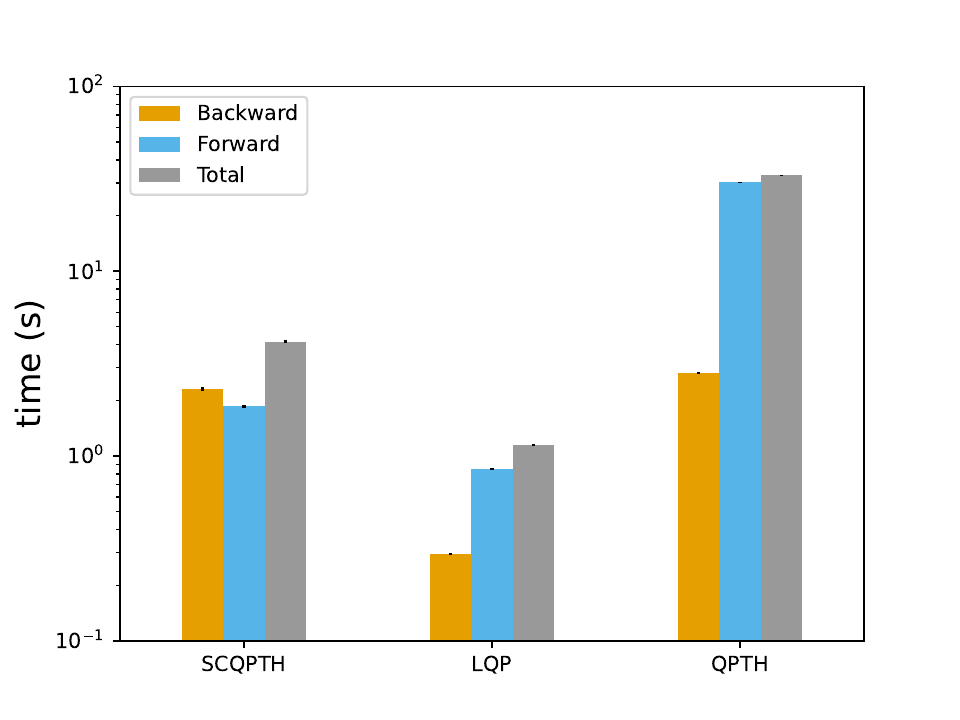}
    \caption{$n = 1000$.}
  \end{subfigure}
  \caption{Computational performance of SCQPTH, LQP, QPTH and Cvxpylayers for box constrained QPs of various problem sizes, $n$, and low stopping tolerance $(1\mathrm{e}{-3})$.}
  \label{fig:exp_1_low}
\end{figure}

\begin{figure}[H]
  \centering
  \begin{subfigure}[b]{0.22\linewidth}
    \includegraphics[width=\linewidth , trim={0mm 0cm 0cm 0cm},clip]{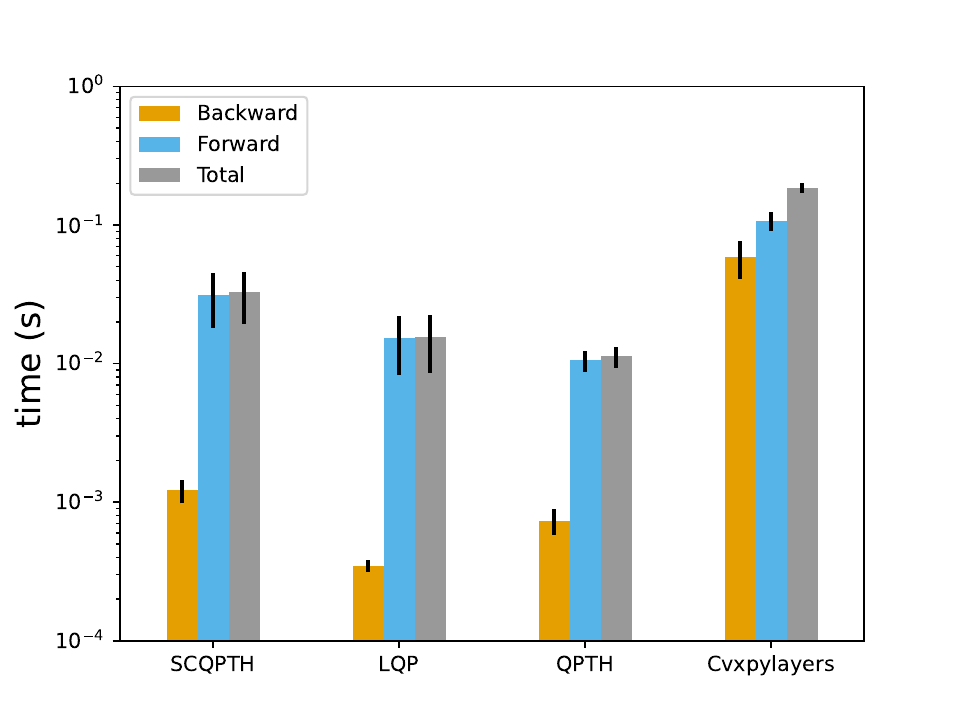}
    \caption{$n = 10$.}
  \end{subfigure}
  \begin{subfigure}[b]{0.22\linewidth}
    \includegraphics[width=\linewidth , trim={0mm 0cm 0cm 0cm},clip]{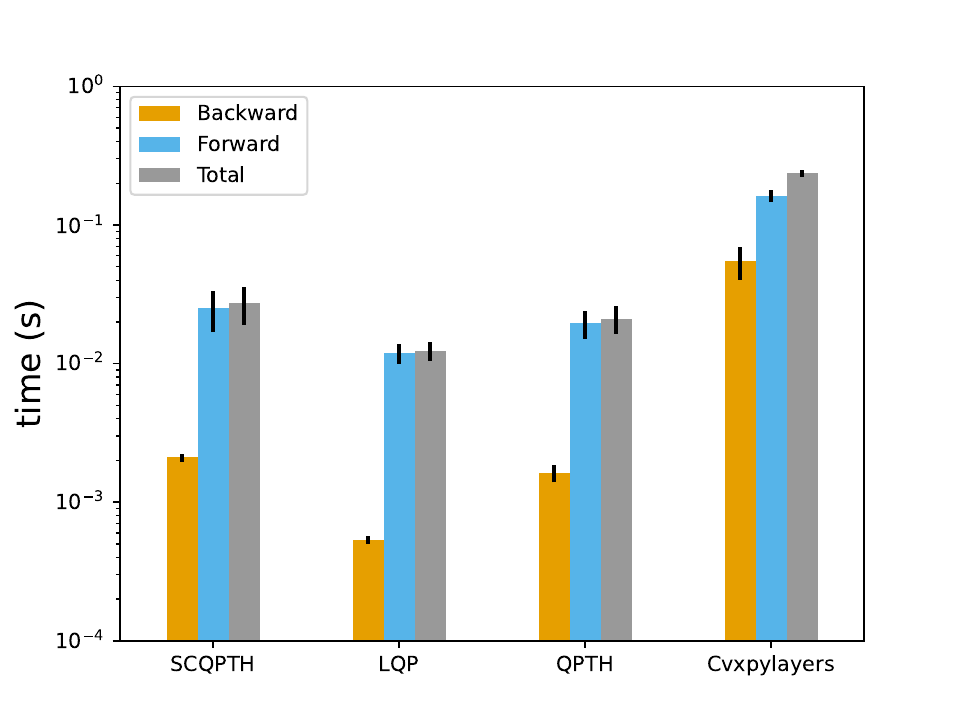}
    \caption{$n = 25$.}
  \end{subfigure}
    \begin{subfigure}[b]{0.22\linewidth}
   \includegraphics[width=\linewidth , trim={0mm 0cm 0cm 0cm},clip]{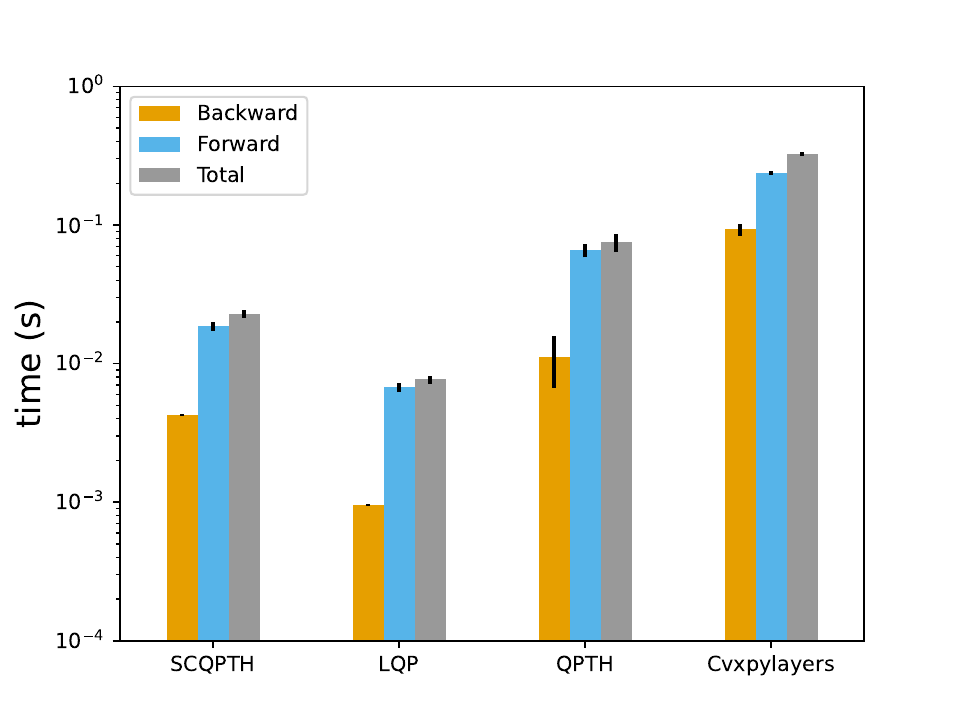}
    \caption{$n= 50$.}
  \end{subfigure}
  \begin{subfigure}[b]{0.22\linewidth}
    \includegraphics[width=\linewidth , trim={0mm 0cm 0cm 0cm},clip]{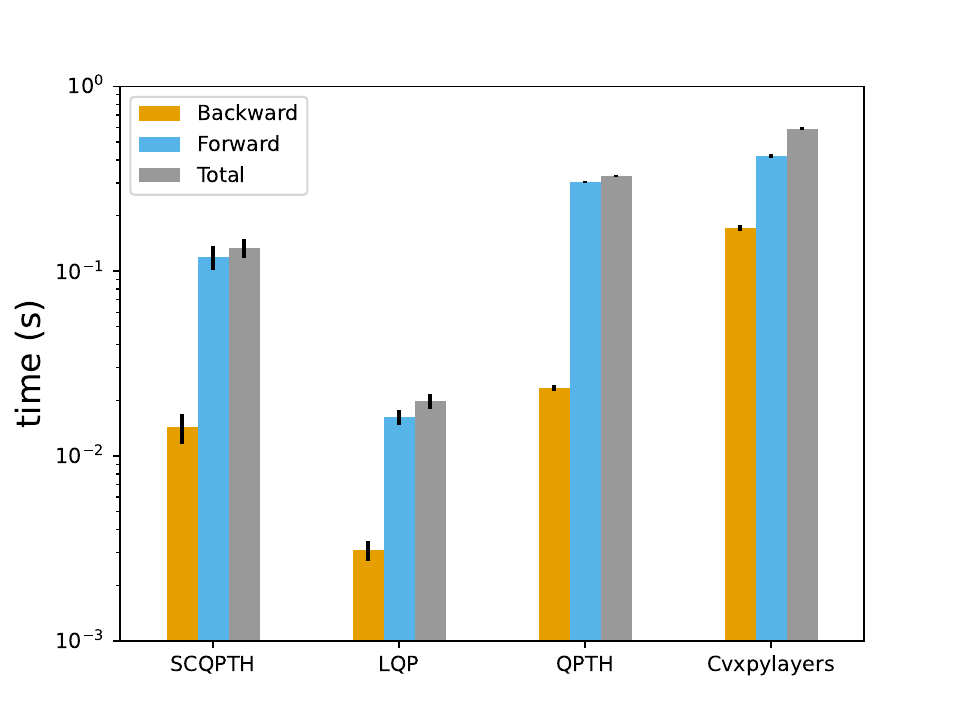}
    \caption{$n = 100$.}
  \end{subfigure}
  \begin{subfigure}[b]{0.22\linewidth}
    \includegraphics[width=\linewidth , trim={0mm 0cm 0cm 0cm},clip]{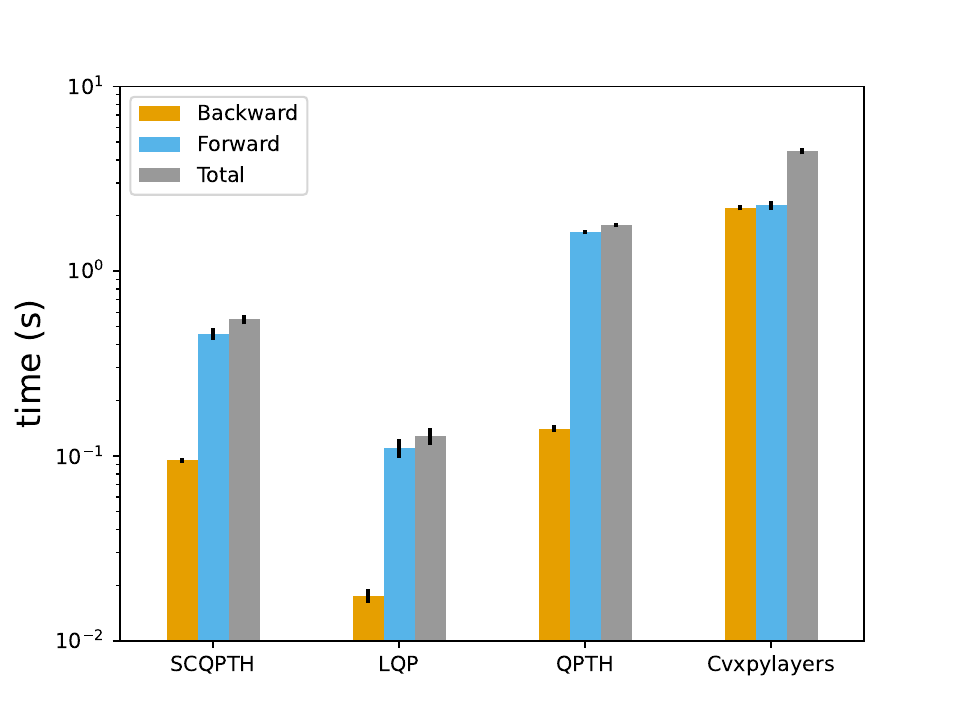}
    \caption{$n = 250$.}
  \end{subfigure}
  \begin{subfigure}[b]{0.22\linewidth}
    \includegraphics[width=\linewidth , trim={0mm 0cm 0cm 0cm},clip]{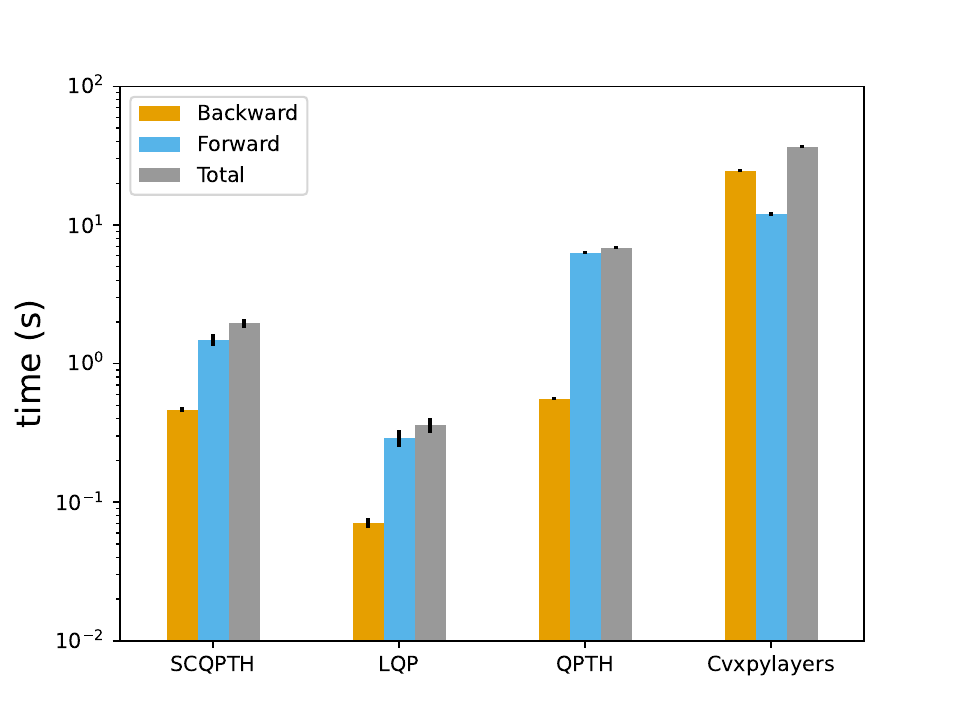}
    \caption{$n = 500$.}
  \end{subfigure}
  \begin{subfigure}[b]{0.22\linewidth}
    \includegraphics[width=\linewidth , trim={0mm 0cm 0cm 0cm},clip]{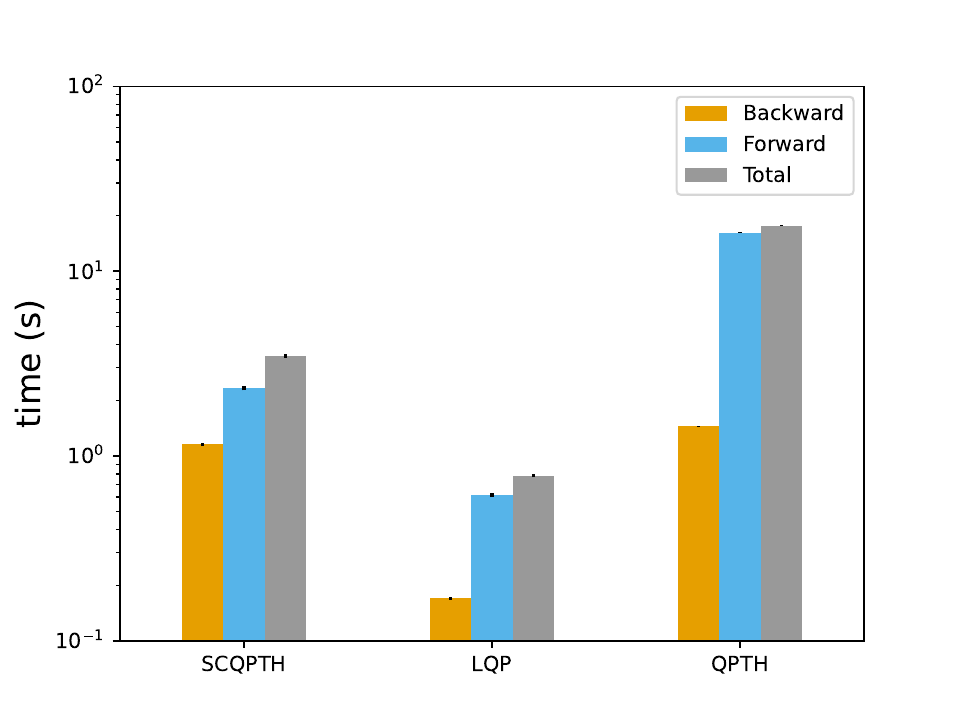}
    \caption{$n = 750$.}
  \end{subfigure}
  \begin{subfigure}[b]{0.22\linewidth}
    \includegraphics[width=\linewidth , trim={0mm 0cm 0cm 0cm},clip]{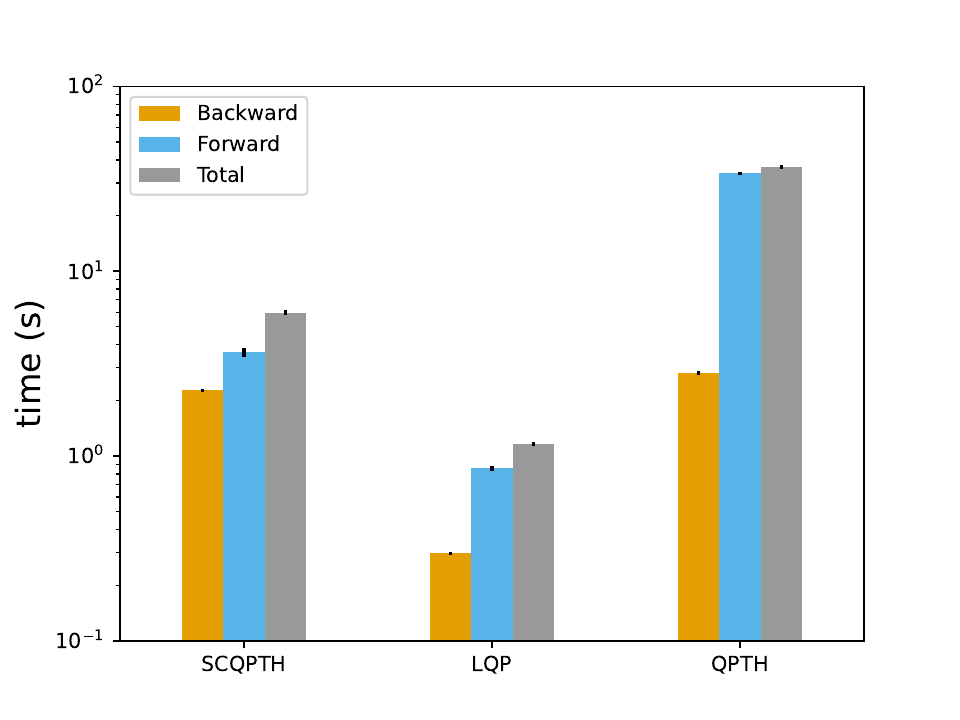}
    \caption{$n = 1000$.}
  \end{subfigure}
  \caption{Computational performance of SCQPTH, LQP, QPTH and Cvxpylayers for box constrained QPs of various problem sizes, $n$, and high stopping tolerance $(1\mathrm{e}{-5})$.}
  \label{fig:exp_1_high}
\end{figure}

Figure \ref{fig:exp_1_low} and Figure \ref{fig:exp_1_high} provide the median runtime and $95\%$-ile confidence interval with a low and high stopping tolerance, respectively. We observe that the LQP layer, which is customized for box constrained QPs, is generally the most computationally efficient method across all problem sizes. The exception is when $n=10$, in which QPTH is computationally most efficient. This is consistent with the observation that interior-point methods can efficiently produce high accurate solutions for problems in low dimensions. In general, the SCQPTH layer is the second most efficient method and provides anywhere from a $1\times$ to over $10\times$ improvement in computational efficiency compared to QPTH and Cvxpylayers. For example, when $n=500$ and stopping tolerance is low, the SCQPTH layer has a median total runtime of $1.35$ seconds whereas QPTH and Cvxpylayers have a median total runtime of $6.40$ seconds and $27.80$ seconds, respectively. The LQP method has a median total runtime of $0.30$ seconds.

\subsection{Experiment 2: Random Constrained QPs}\label{sec:results_2}
We generate randomly constrained QPs of dimension: $n \in \{10, 25, 50, 100, 250, 500, 740, 1000\}$ and consider low $(1\mathrm{e}{-3})$ and high $(1\mathrm{e}{-5})$ absolute and relative stopping tolerances. Problem variables are generated as follows. As before, we set $\bQ =  \bL^T\bL + 0.01 \bI$ where $\bL \in \mathbb{R}^{n \times n}$ and entries $\bL_{ij} \sim \mathcal{N}(0, 1)$ with $50\%$ probability of being non-zero,  $\bp_i \sim \mathcal{N}(0, 1)$, $\bl_j \sim \mathcal{U}(-1,0)$ and $\bu_j \sim \mathcal{U}(0,1)$. We randomly generate $\bA \in \mathbb{R}^{m \times n}$ with entries $\bA_{ij} \sim \mathcal{N}(0, 1)$ with $15\%$ probability of being non-zero and consider the case where $m = n$ and $m = 2n$. Experiment results are averaged over 10 independent trials with a batch size of 32.

\begin{figure}[H]
  \centering
  \begin{subfigure}[b]{0.22\linewidth}
    \includegraphics[width=\linewidth , trim={0mm 0cm 0cm 0cm},clip]{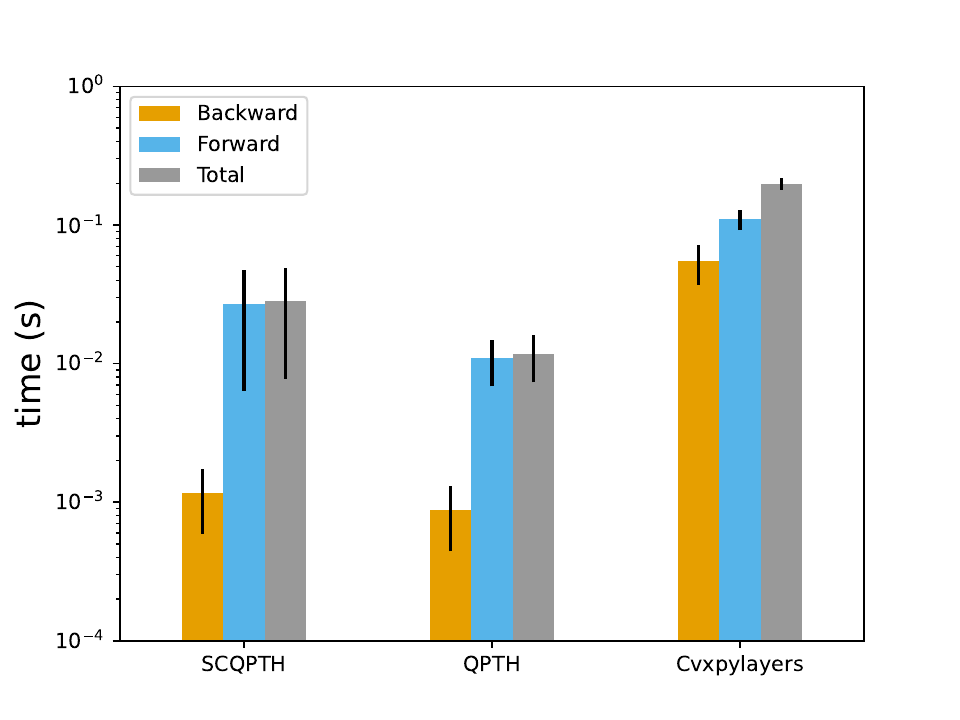}
    \caption{$n = 10, m=10$.}
  \end{subfigure}
  \begin{subfigure}[b]{0.22\linewidth}
    \includegraphics[width=\linewidth , trim={0mm 0cm 0cm 0cm},clip]{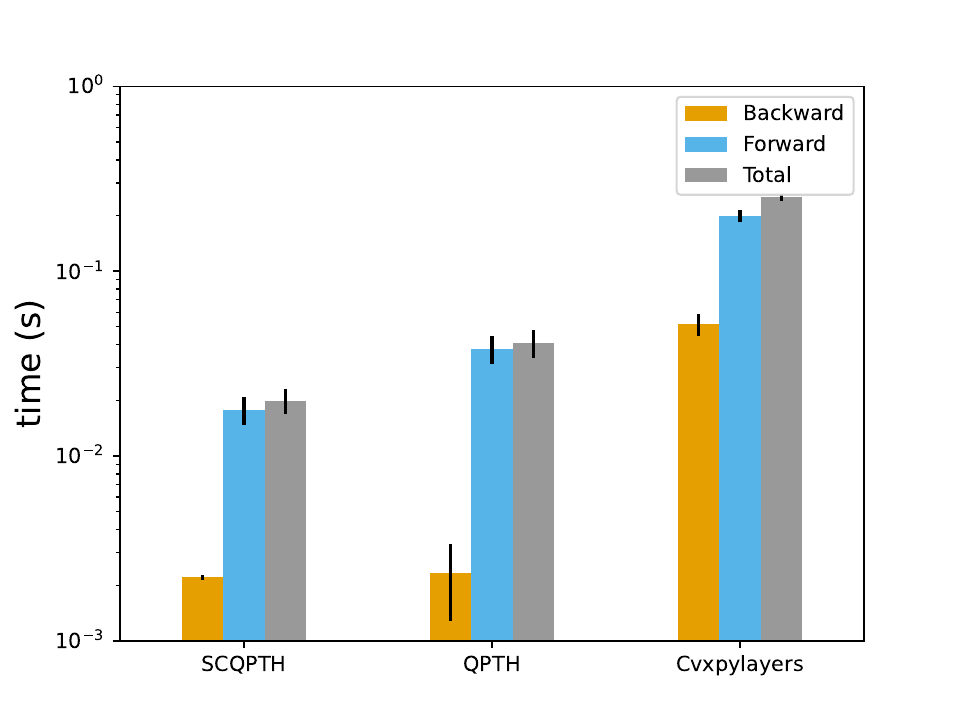}
    \caption{$n = 25, m=25$.}
  \end{subfigure}
    \begin{subfigure}[b]{0.22\linewidth}
   \includegraphics[width=\linewidth , trim={0mm 0cm 0cm 0cm},clip]{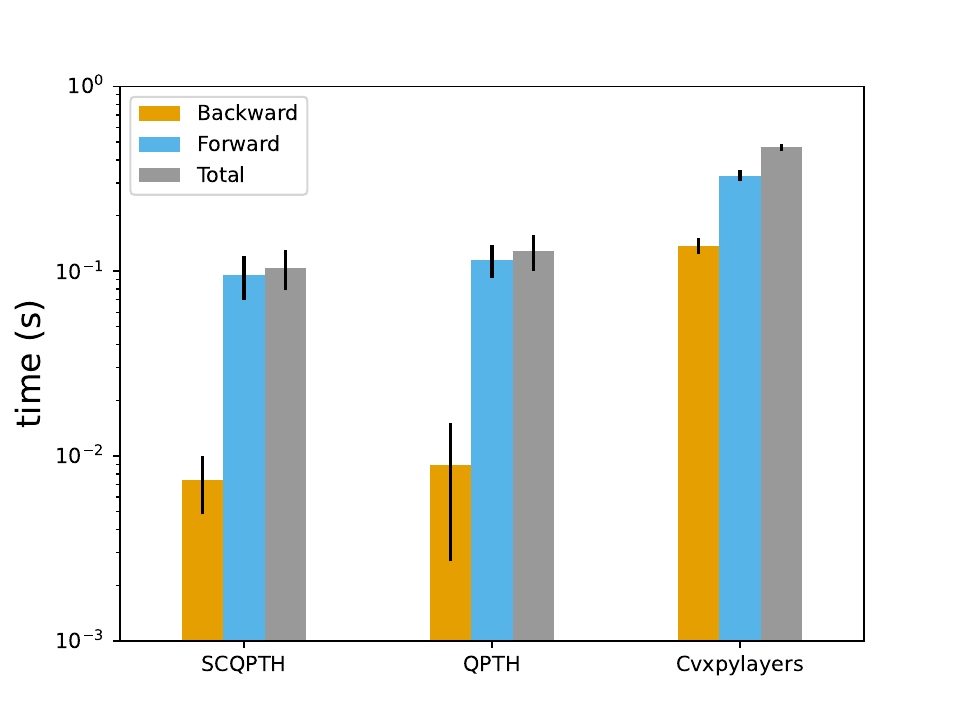}
    \caption{$n= 50, m=50$.}
  \end{subfigure}
  \begin{subfigure}[b]{0.22\linewidth}
    \includegraphics[width=\linewidth , trim={0mm 0cm 0cm 0cm},clip]{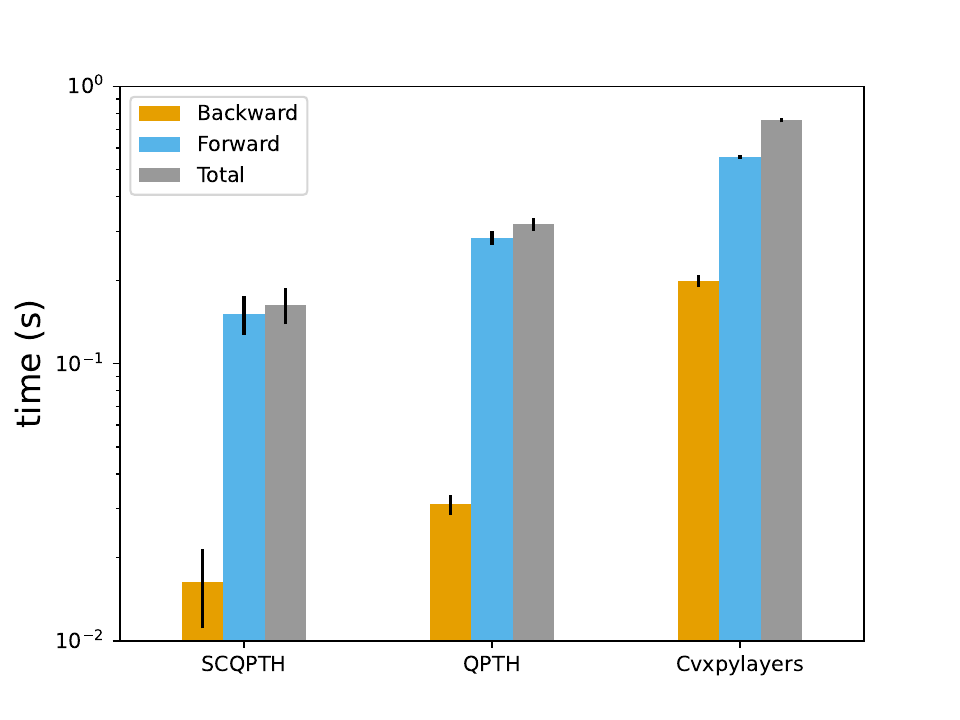}
    \caption{$n = 100, m=100$.}
  \end{subfigure}
  \begin{subfigure}[b]{0.22\linewidth}
    \includegraphics[width=\linewidth , trim={0mm 0cm 0cm 0cm},clip]{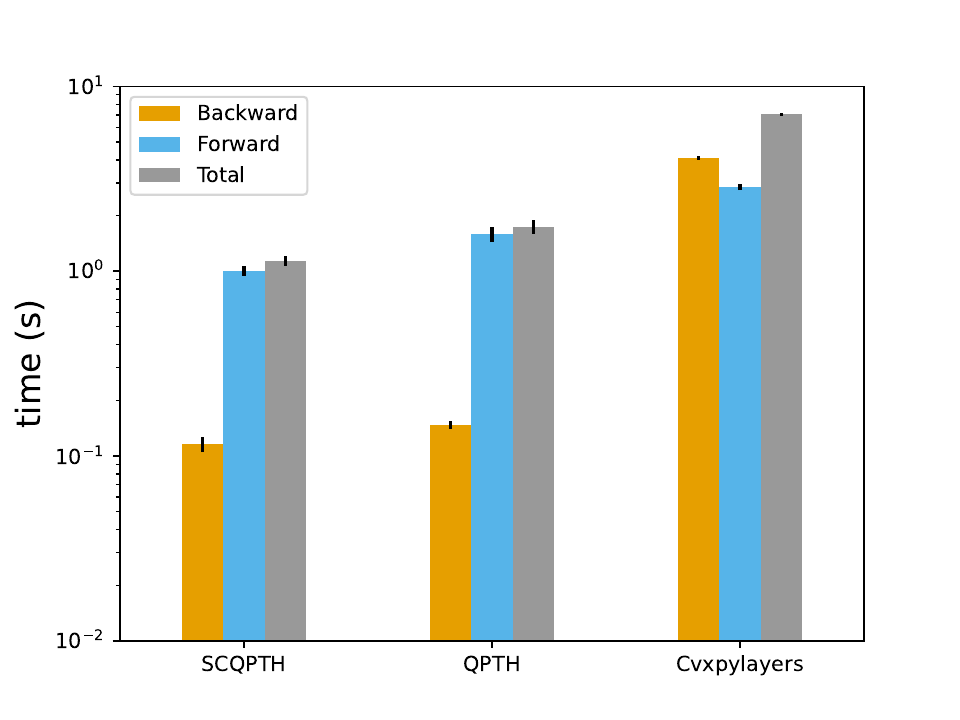}
    \caption{$n = 250, m=250$.}
  \end{subfigure}
  \begin{subfigure}[b]{0.22\linewidth}
    \includegraphics[width=\linewidth , trim={0mm 0cm 0cm 0cm},clip]{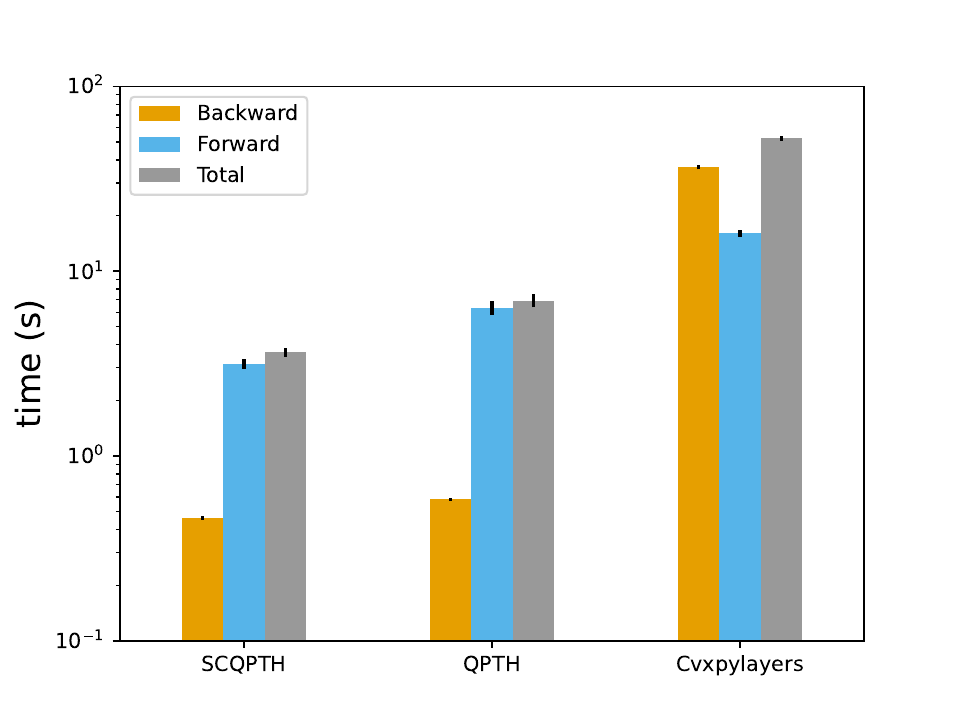}
    \caption{$n = 500, m=500$.}
  \end{subfigure}
  \begin{subfigure}[b]{0.22\linewidth}
    \includegraphics[width=\linewidth , trim={0mm 0cm 0cm 0cm},clip]{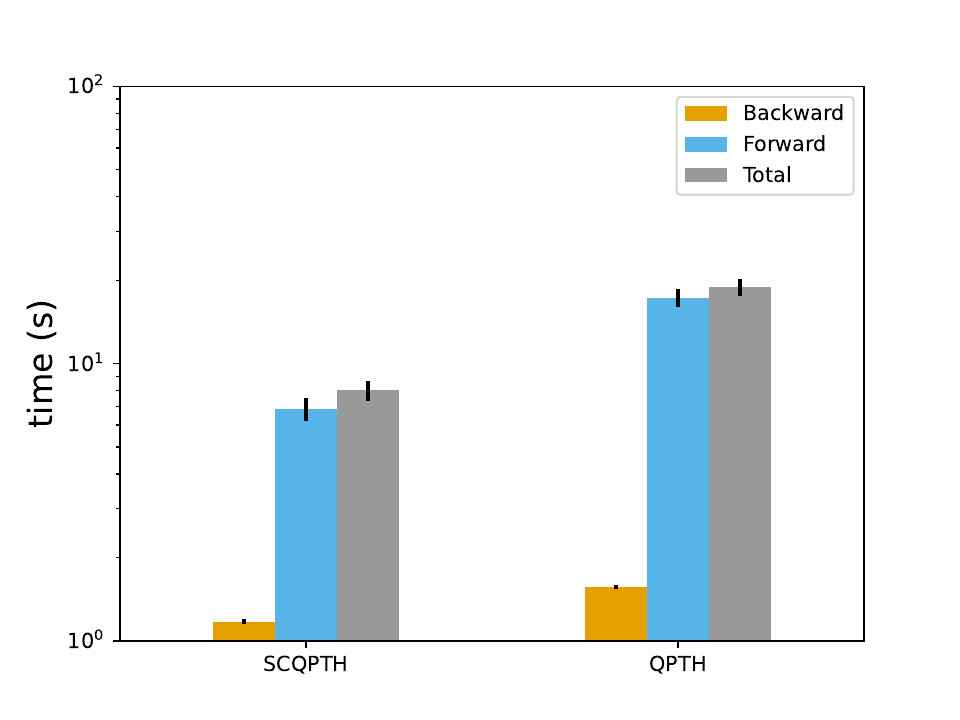}
    \caption{$n = 750, m=750$.}
  \end{subfigure}
  \begin{subfigure}[b]{0.22\linewidth}
    \includegraphics[width=\linewidth , trim={0mm 0cm 0cm 0cm},clip]{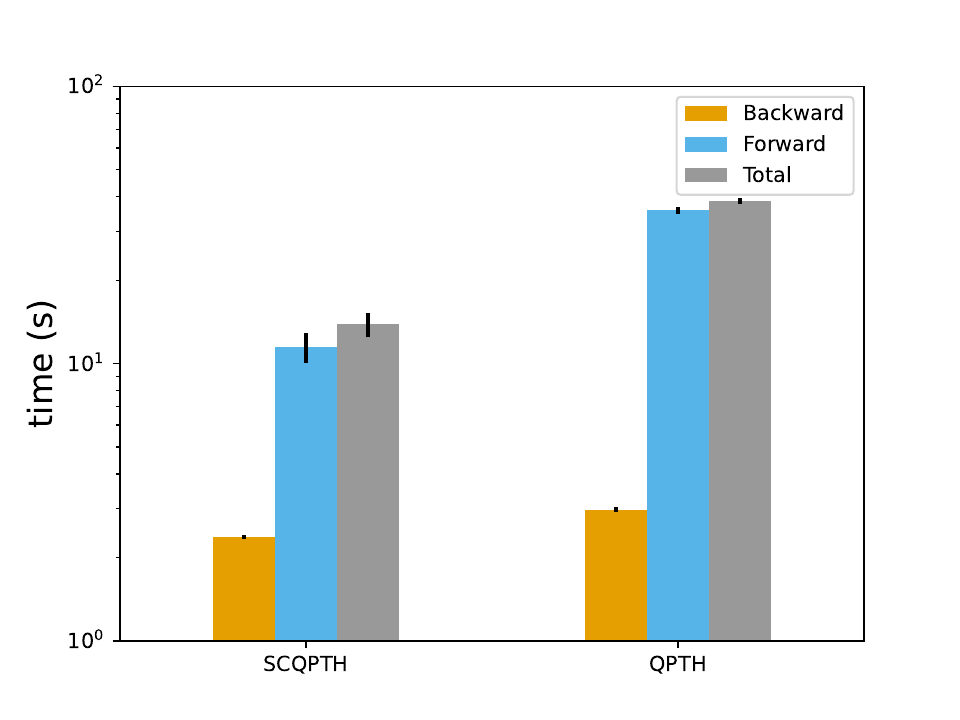}
    \caption{$n = 1000, m=1000$.}
  \end{subfigure}
  \caption{Computational performance of SCQPTH,  QPTH and Cvxpylayers for random QPs of various problem sizes, $n$, constraints $m=n$, and low stopping tolerance $(1\mathrm{e}{-3})$.}
  \label{fig:exp_2_low}
\end{figure}

Figure \ref{fig:exp_2_low} provides the median runtime and $95\%$-ile confidence interval with a low stopping tolerance and number of constraints $m=n$. As before the SCQPTH layer is the most efficient method and for large scale QPs provides a $1\times$ to over $10\times$ improvement in computational efficiency compared to QPTH and Cvxpylayers.  For example, when $n=500$ the SCQPTH layer has a median total runtime of $3.65$ seconds whereas QPTH and Cvxpylayers have a median total runtime of $6.95$ seconds and $52.40$ seconds, respectively.

\begin{figure}[H]
  \centering
  \begin{subfigure}[b]{0.22\linewidth}
    \includegraphics[width=\linewidth , trim={0mm 0cm 0cm 0cm},clip]{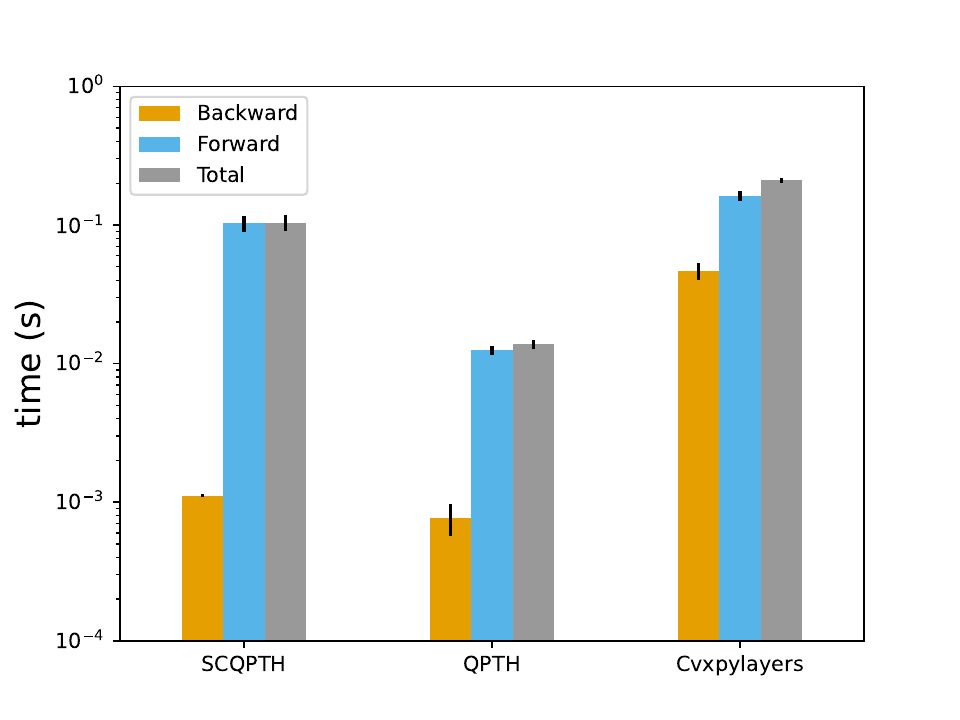}
    \caption{$n = 10, m=10$.}
  \end{subfigure}
  \begin{subfigure}[b]{0.22\linewidth}
    \includegraphics[width=\linewidth , trim={0mm 0cm 0cm 0cm},clip]{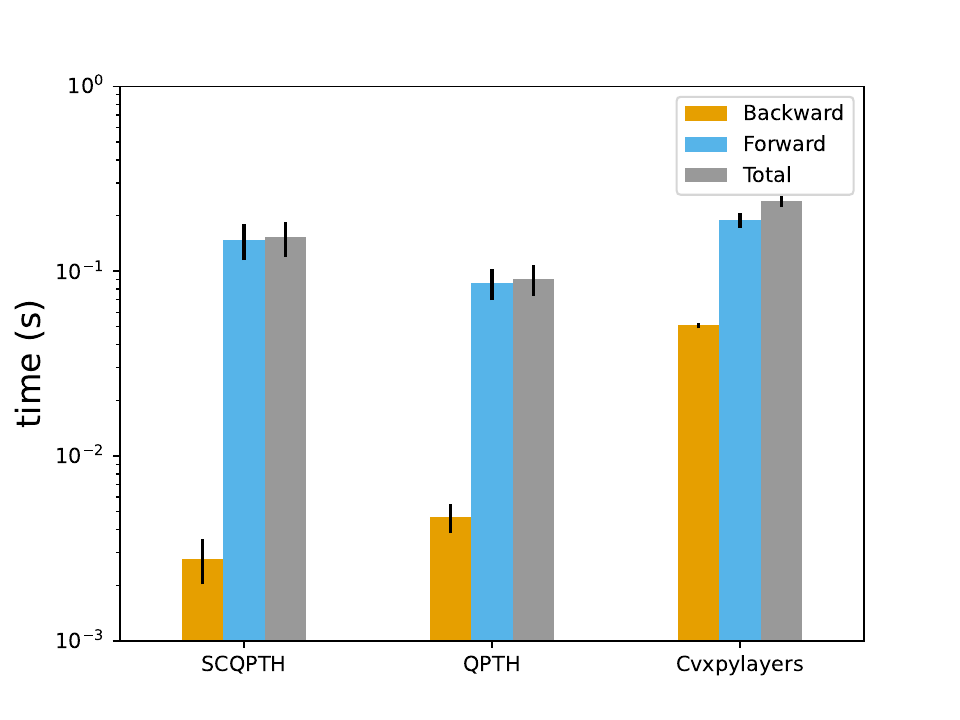}
    \caption{$n = 25, m=25$.}
  \end{subfigure}
    \begin{subfigure}[b]{0.22\linewidth}
   \includegraphics[width=\linewidth , trim={0mm 0cm 0cm 0cm},clip]{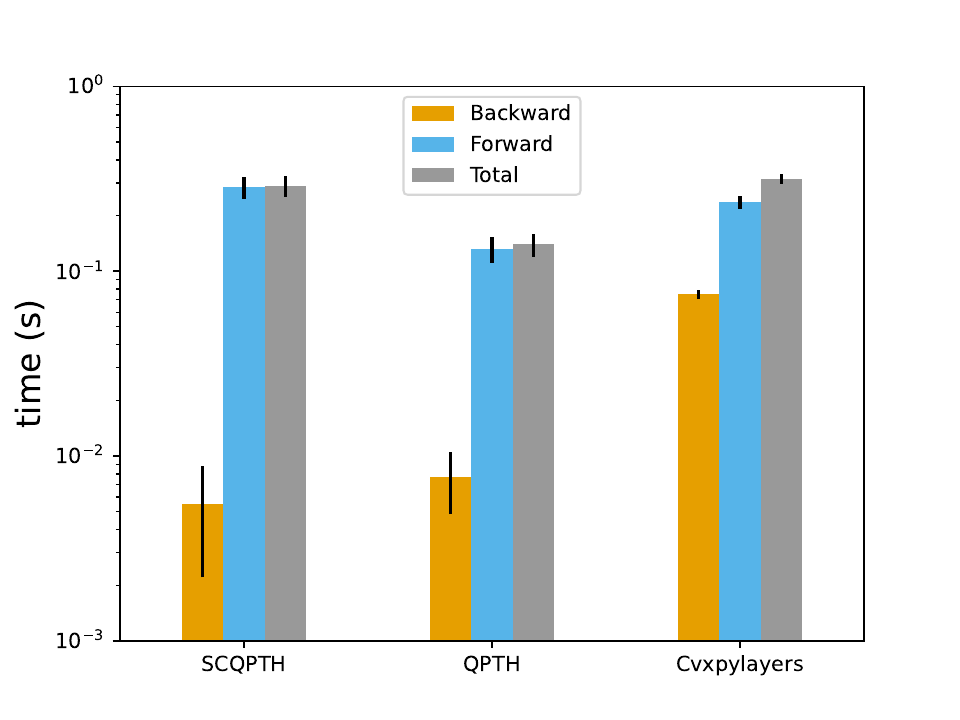}
    \caption{$n= 50, m=50$.}
  \end{subfigure}
  \begin{subfigure}[b]{0.22\linewidth}
    \includegraphics[width=\linewidth , trim={0mm 0cm 0cm 0cm},clip]{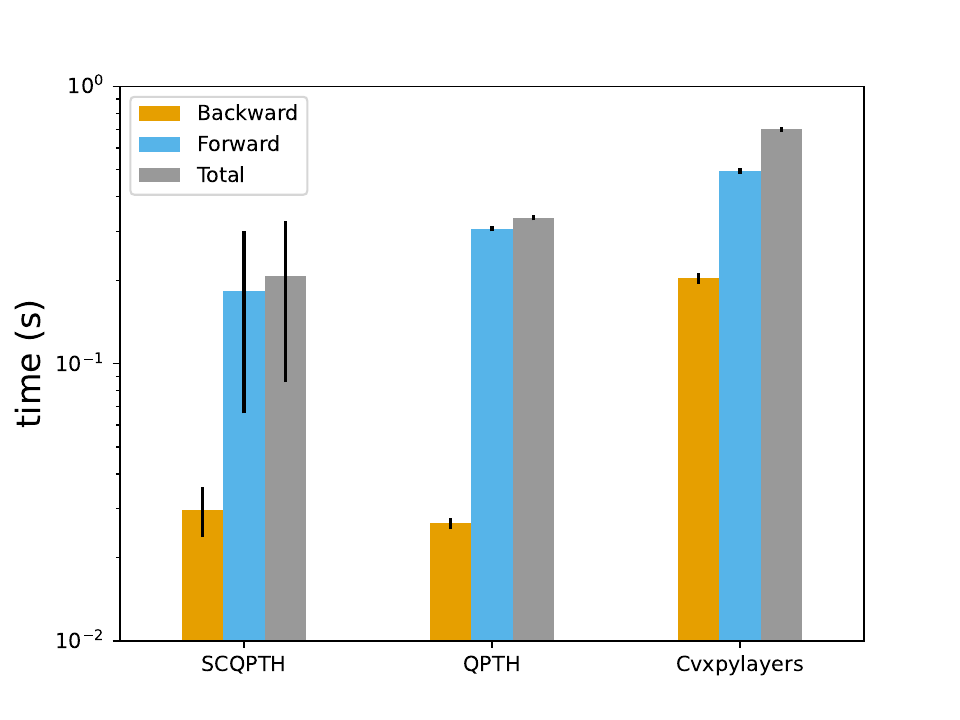}
    \caption{$n = 100, m=100$.}
  \end{subfigure}
  \begin{subfigure}[b]{0.22\linewidth}
    \includegraphics[width=\linewidth , trim={0mm 0cm 0cm 0cm},clip]{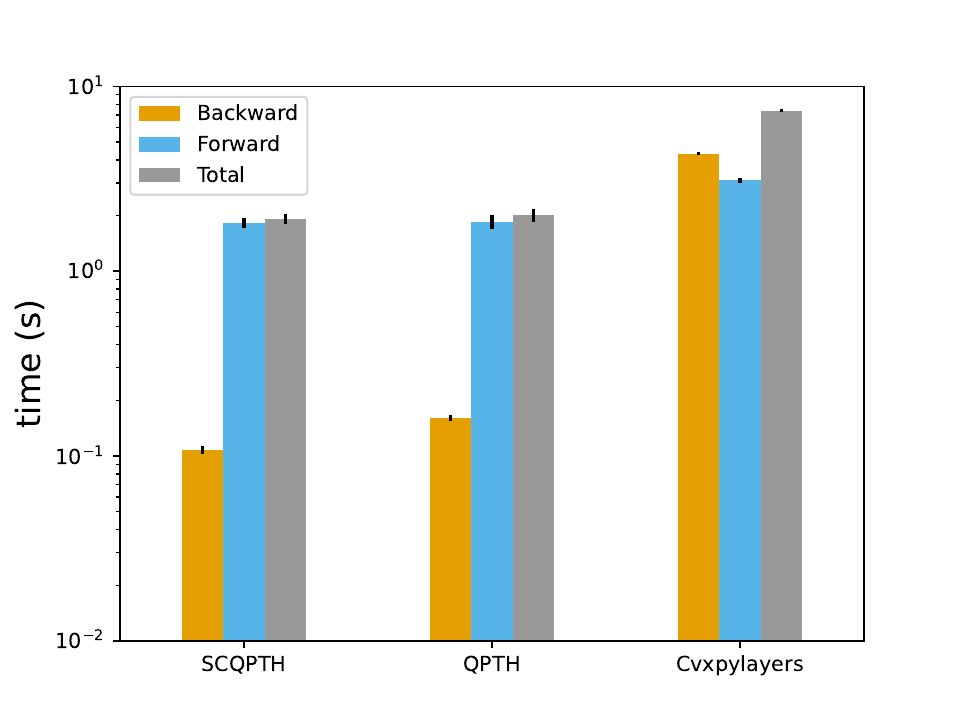}
    \caption{$n = 250,m=250$.}
  \end{subfigure}
  \begin{subfigure}[b]{0.22\linewidth}
    \includegraphics[width=\linewidth , trim={0mm 0cm 0cm 0cm},clip]{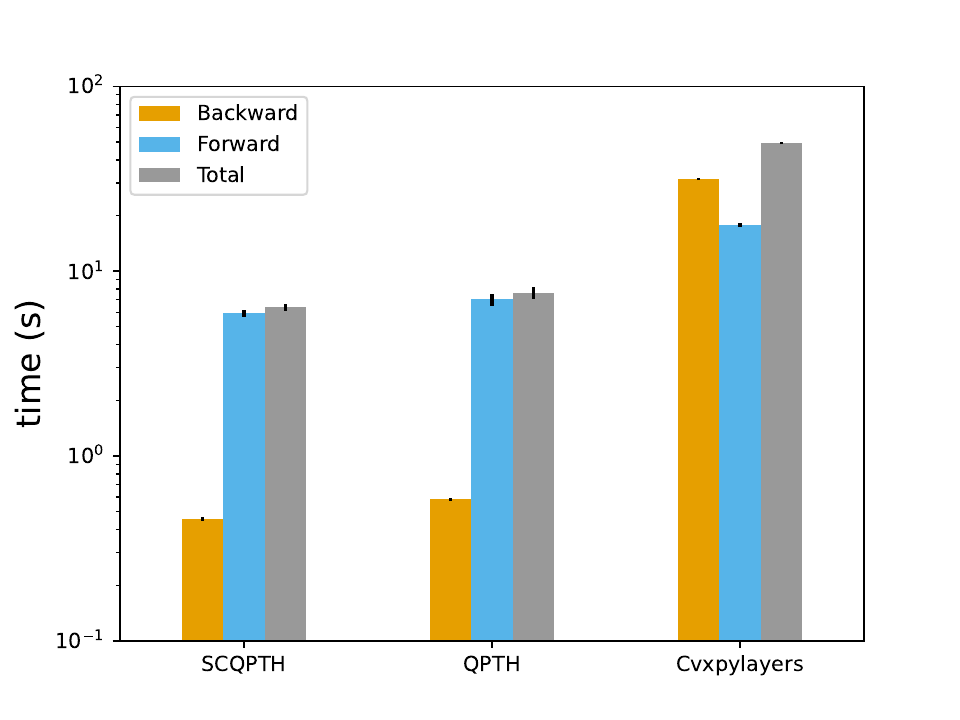}
    \caption{$n = 500,m=500$.}
  \end{subfigure}
  \begin{subfigure}[b]{0.22\linewidth}
    \includegraphics[width=\linewidth , trim={0mm 0cm 0cm 0cm},clip]{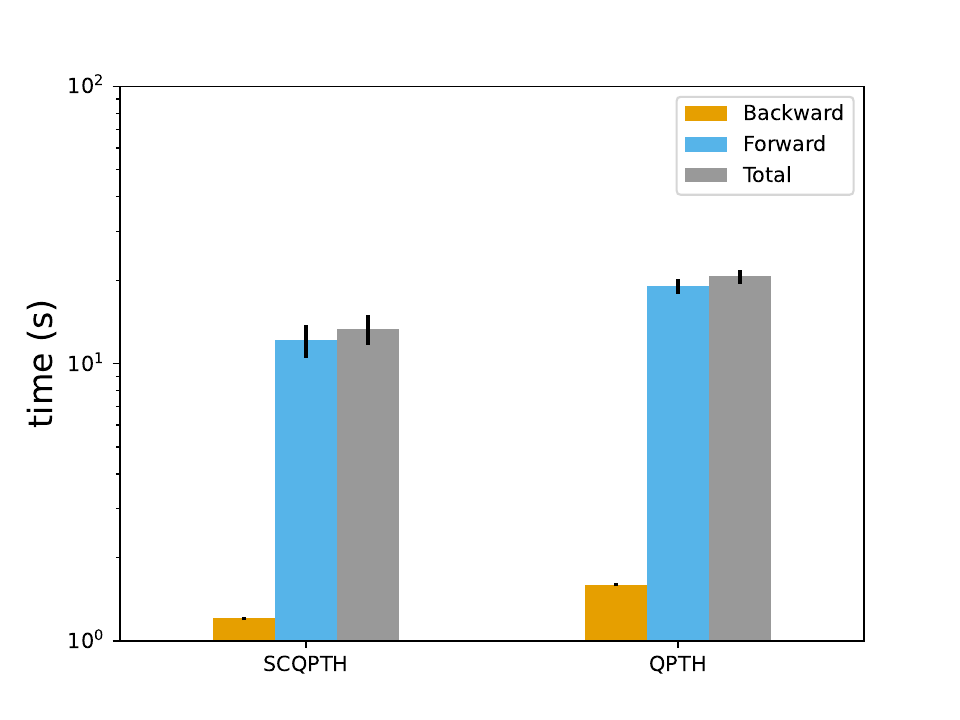}
    \caption{$n = 750, m=750$.}
  \end{subfigure}
  \begin{subfigure}[b]{0.22\linewidth}
    \includegraphics[width=\linewidth , trim={0mm 0cm 0cm 0cm},clip]{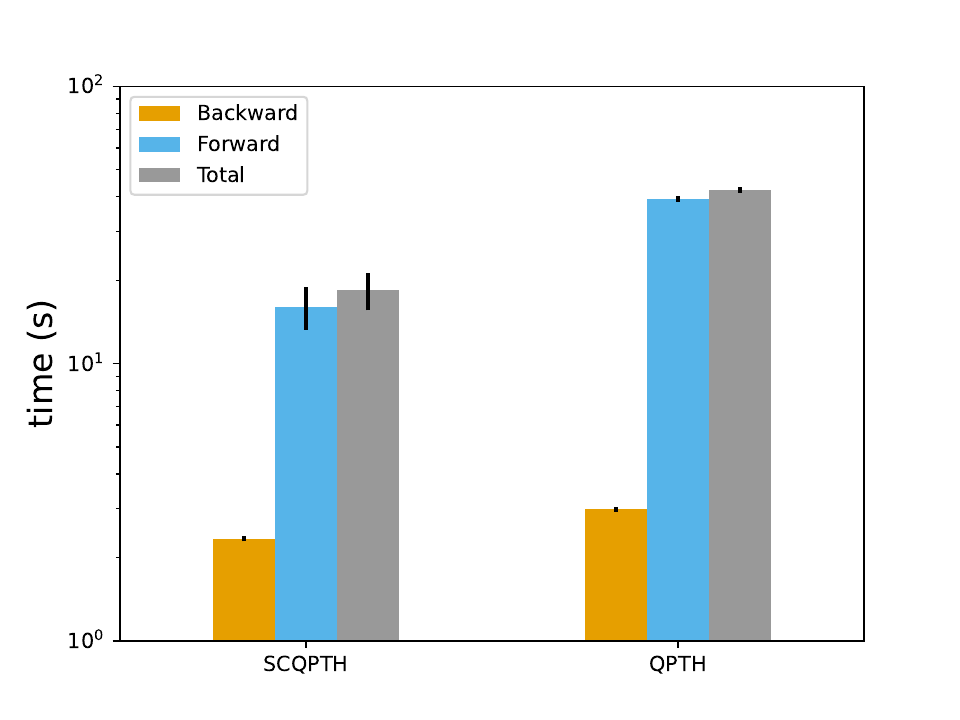}
    \caption{$n = 1000, m=1000$.}
  \end{subfigure}
  \caption{Computational performance of SCQPTH, QPTH and Cvxpylayers for random QPs of various problem sizes, $n$, constraints $m=n$, and high stopping tolerance $(1\mathrm{e}{-5})$.}
  \label{fig:exp_2_high}
\end{figure}

Similarly, Figure \ref{fig:exp_2_high} provides the median runtime and $95\%$-ile confidence interval with a high stopping tolerance and number of constraints $m=n$. In general, for large scale problems, the SCQPTH layer is the most efficient method and provides a more modest $1\times - 7\times$ improvement in computational efficiency compared to QPTH and Cvxpylayers.  For example, when $n=500$ the SCQPTH layer has a median total runtime of $6.40$ seconds whereas QPTH and Cvxpylayers have a median total runtime of $7.60$ seconds and $49.30$ seconds, respectively.

\begin{figure}[H]
  \centering
  \begin{subfigure}[b]{0.22\linewidth}
    \includegraphics[width=\linewidth , trim={0mm 0cm 0cm 0cm},clip]{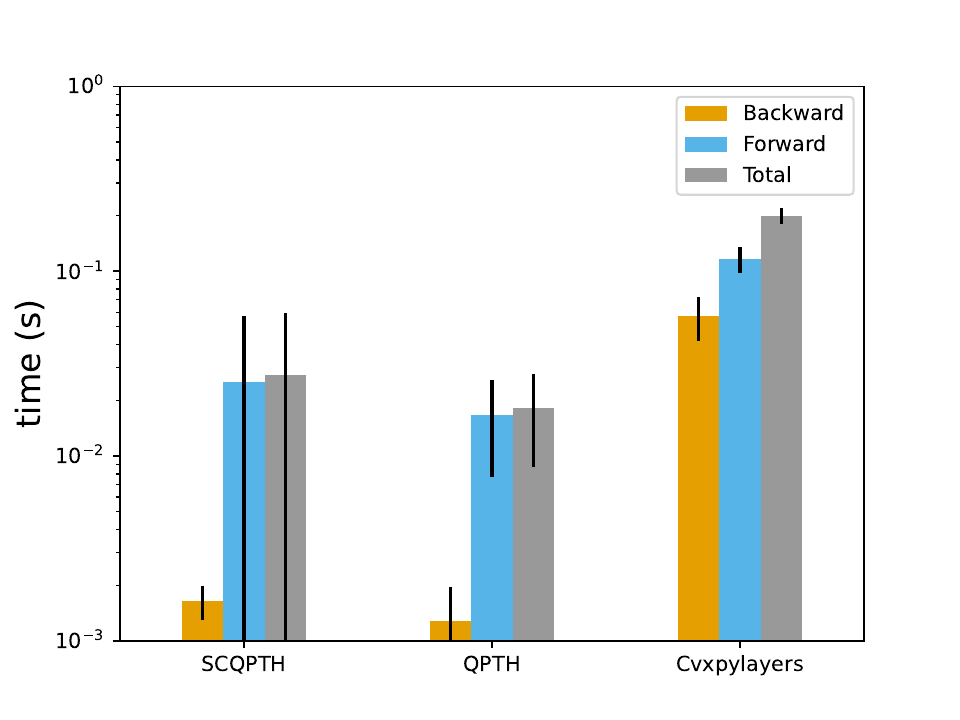}
    \caption{$n = 10, m = 20$.}
  \end{subfigure}
  \begin{subfigure}[b]{0.22\linewidth}
    \includegraphics[width=\linewidth , trim={0mm 0cm 0cm 0cm},clip]{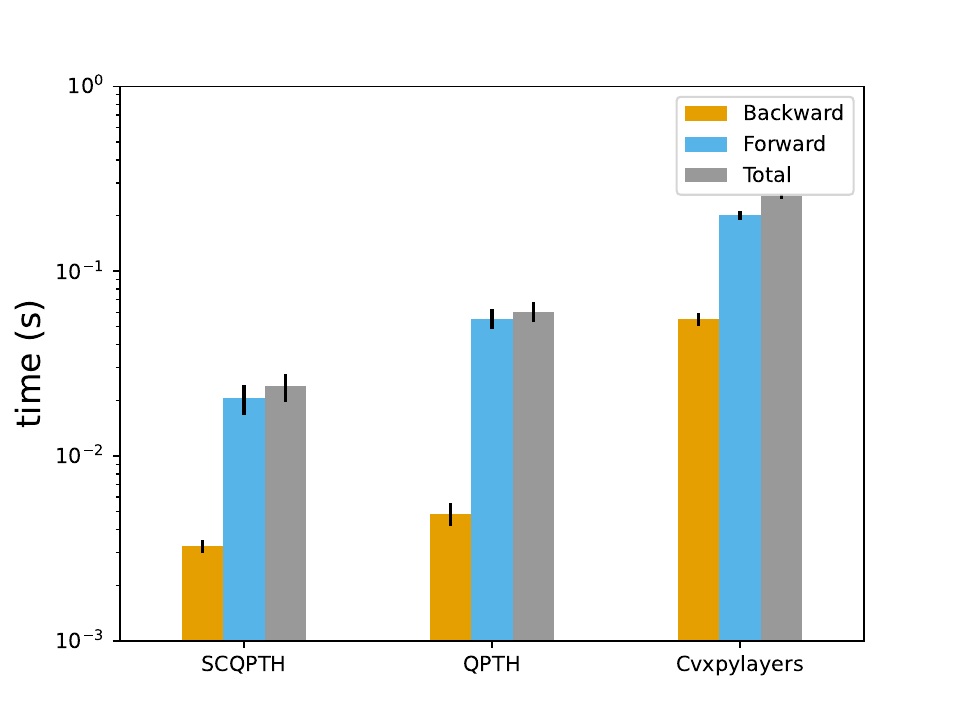}
    \caption{$n = 25, m=50$.}
  \end{subfigure}
    \begin{subfigure}[b]{0.22\linewidth}
   \includegraphics[width=\linewidth , trim={0mm 0cm 0cm 0cm},clip]{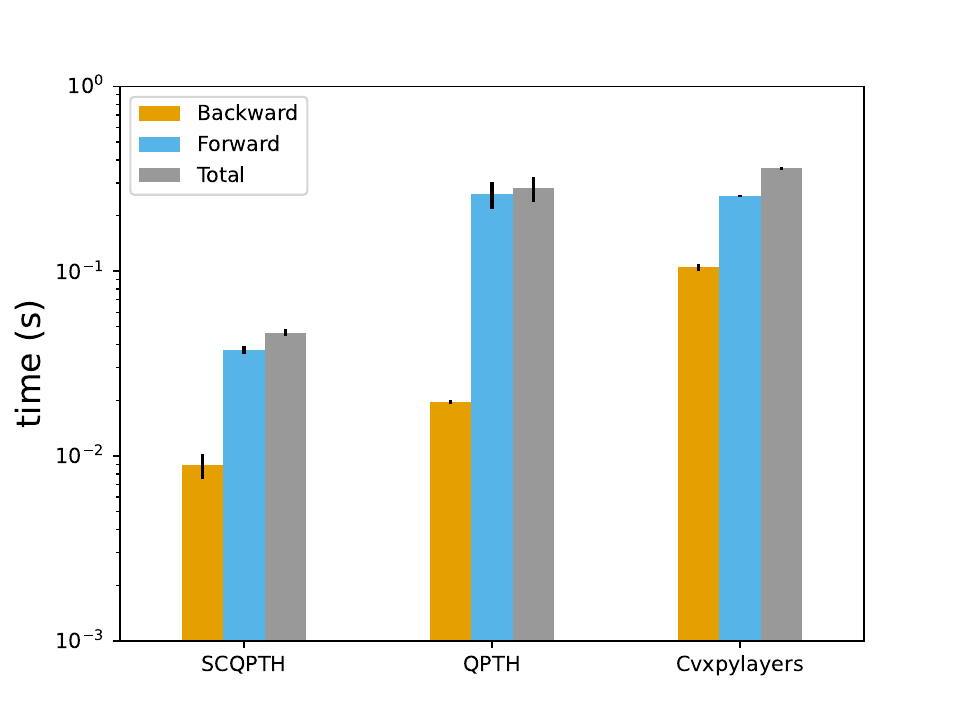}
    \caption{$n= 50, m=100$.}
  \end{subfigure}
  \begin{subfigure}[b]{0.22\linewidth}
    \includegraphics[width=\linewidth , trim={0mm 0cm 0cm 0cm},clip]{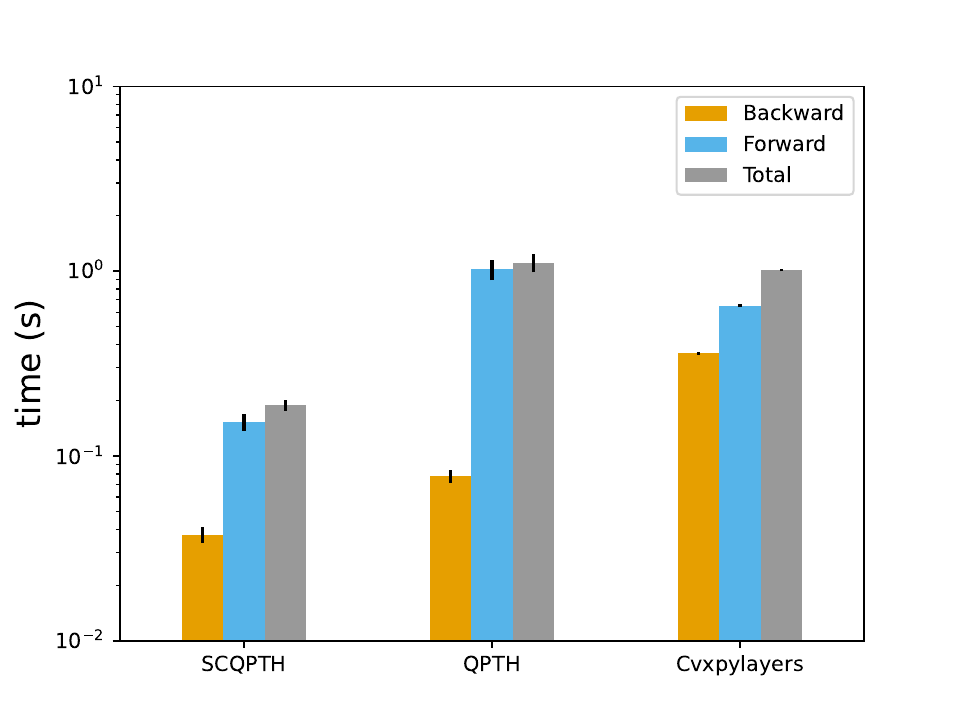}
    \caption{$n = 100, m=200$.}
  \end{subfigure}
  \begin{subfigure}[b]{0.22\linewidth}
    \includegraphics[width=\linewidth , trim={0mm 0cm 0cm 0cm},clip]{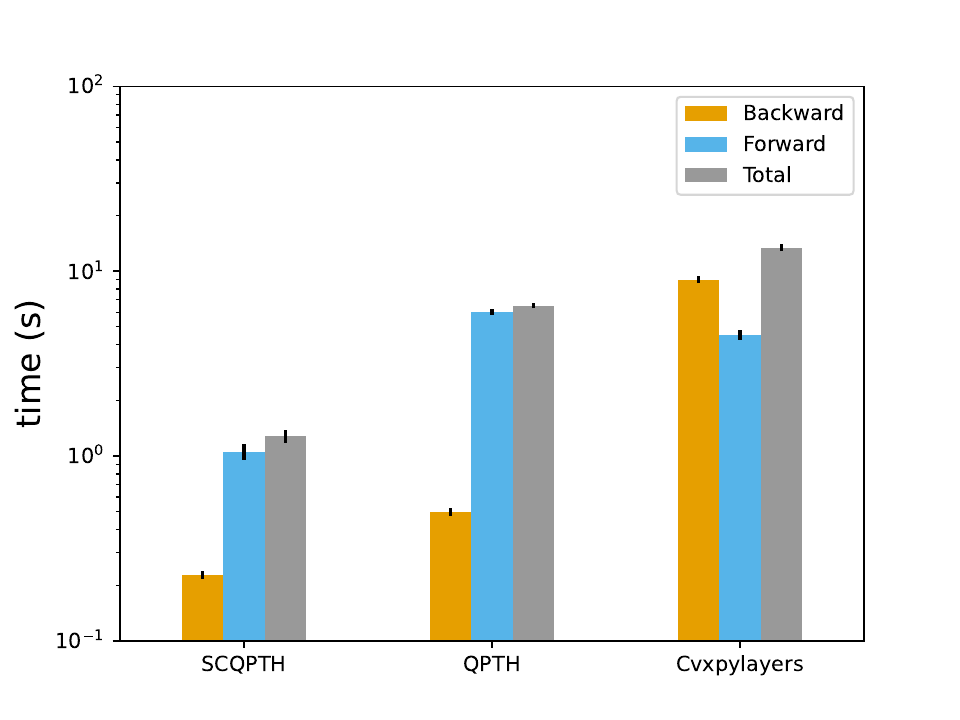}
    \caption{$n = 250, m=500$.}
  \end{subfigure}
  \begin{subfigure}[b]{0.22\linewidth}
    \includegraphics[width=\linewidth , trim={0mm 0cm 0cm 0cm},clip]{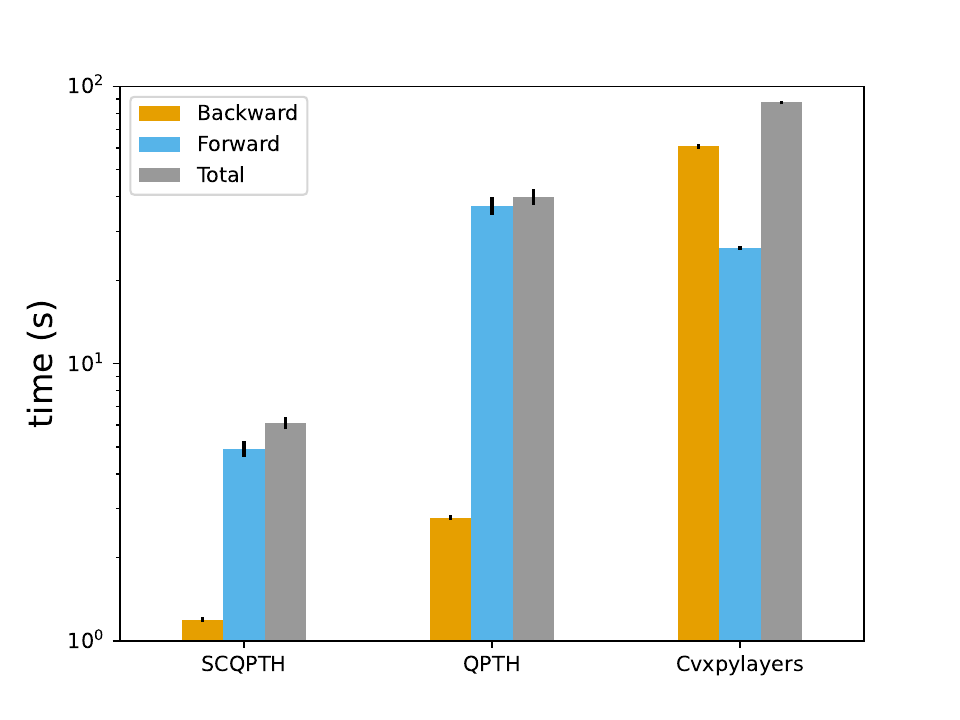}
    \caption{$n = 500, m=1000$.}
  \end{subfigure}
  \begin{subfigure}[b]{0.22\linewidth}
    \includegraphics[width=\linewidth , trim={0mm 0cm 0cm 0cm},clip]{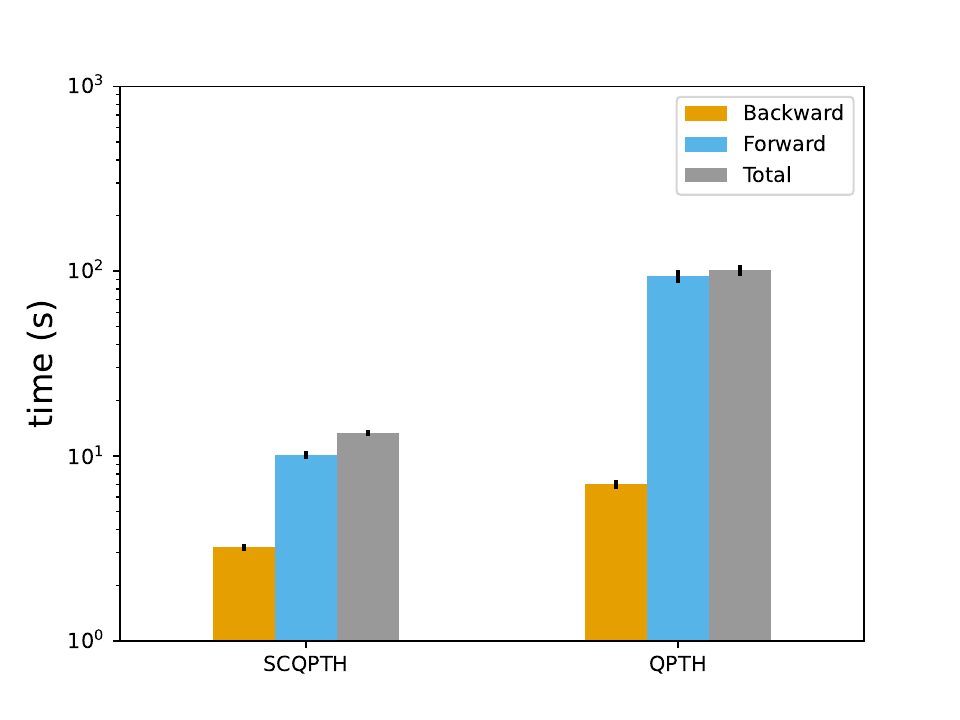}
    \caption{$n = 750, m=1500$.}
  \end{subfigure}
  \begin{subfigure}[b]{0.22\linewidth}
    \includegraphics[width=\linewidth , trim={0mm 0cm 0cm 0cm},clip]{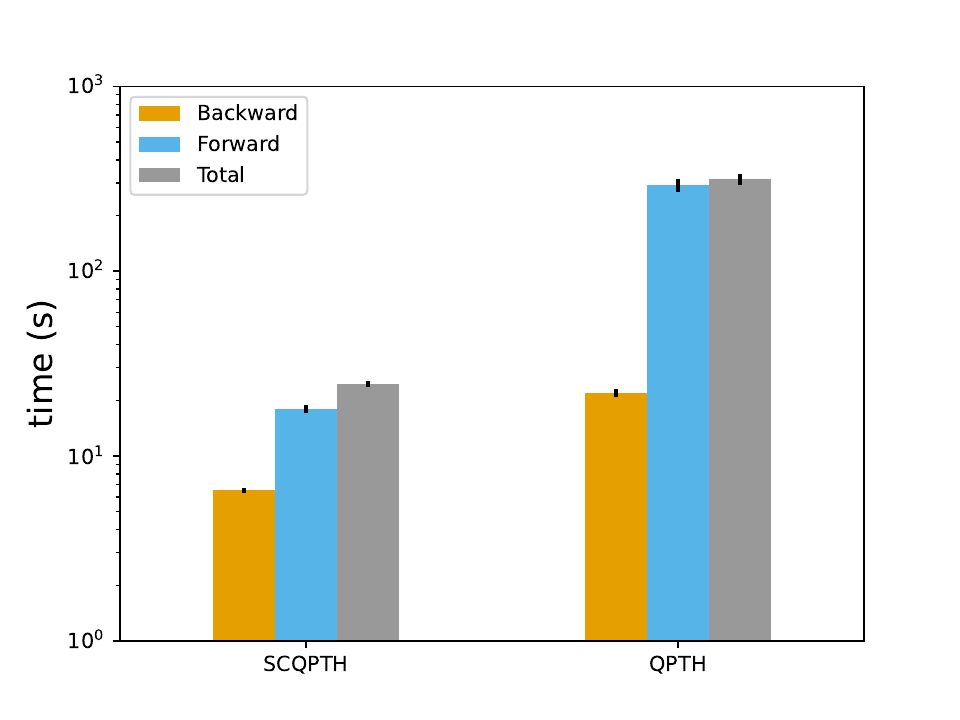}
    \caption{$n = 1000, m=2000$.}
  \end{subfigure}
  \caption{Computational performance of SCQPTH, QPTH and Cvxpylayers for random QPs of various problem sizes, $n$, constraints $m=2n$, and low stopping tolerance $(1\mathrm{e}{-3})$.}
  \label{fig:exp_2_low_2}
\end{figure}

\begin{figure}[h]
  \centering
  \begin{subfigure}[b]{0.22\linewidth}
    \includegraphics[width=\linewidth , trim={0mm 0cm 0cm 0cm},clip]{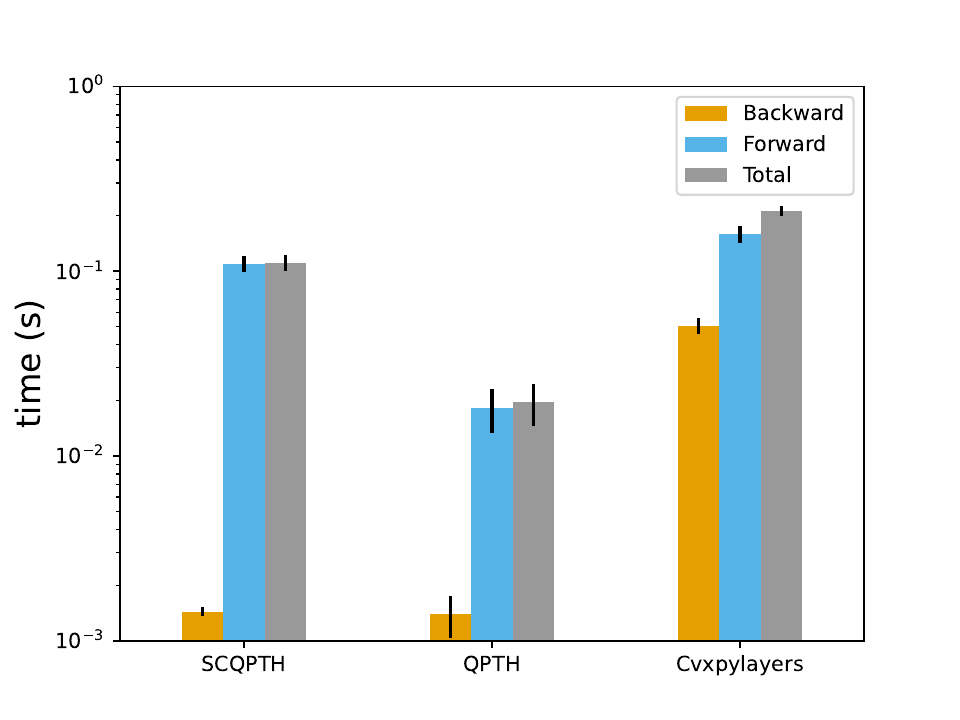}
    \caption{$n = 10, m=20$.}
  \end{subfigure}
  \begin{subfigure}[b]{0.22\linewidth}
    \includegraphics[width=\linewidth , trim={0mm 0cm 0cm 0cm},clip]{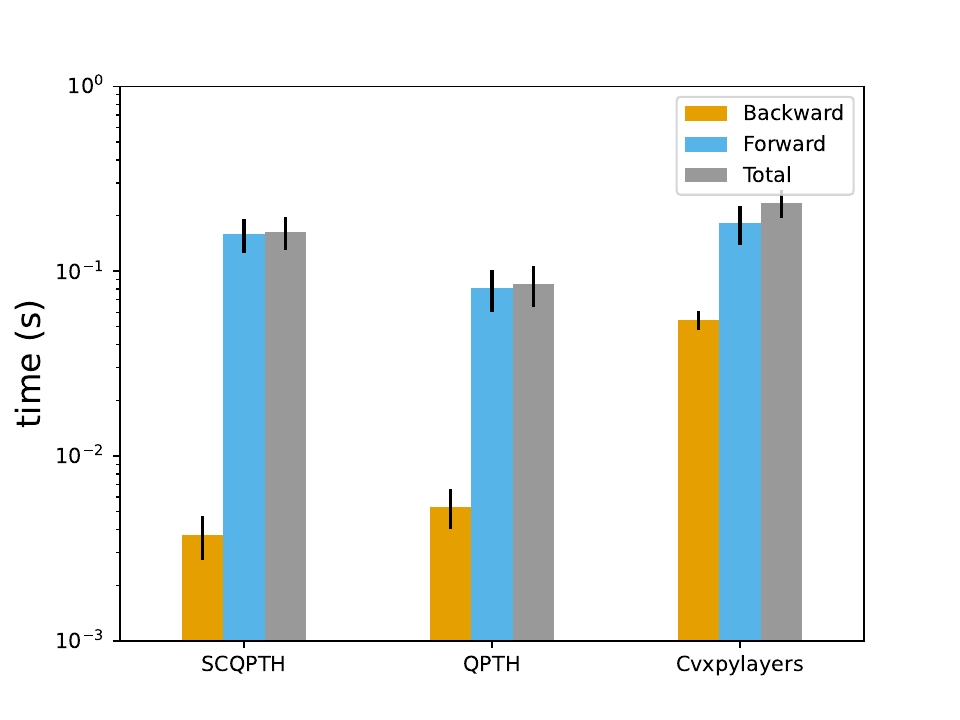}
    \caption{$n = 25, m=50$.}
  \end{subfigure}
    \begin{subfigure}[b]{0.22\linewidth}
   \includegraphics[width=\linewidth , trim={0mm 0cm 0cm 0cm},clip]{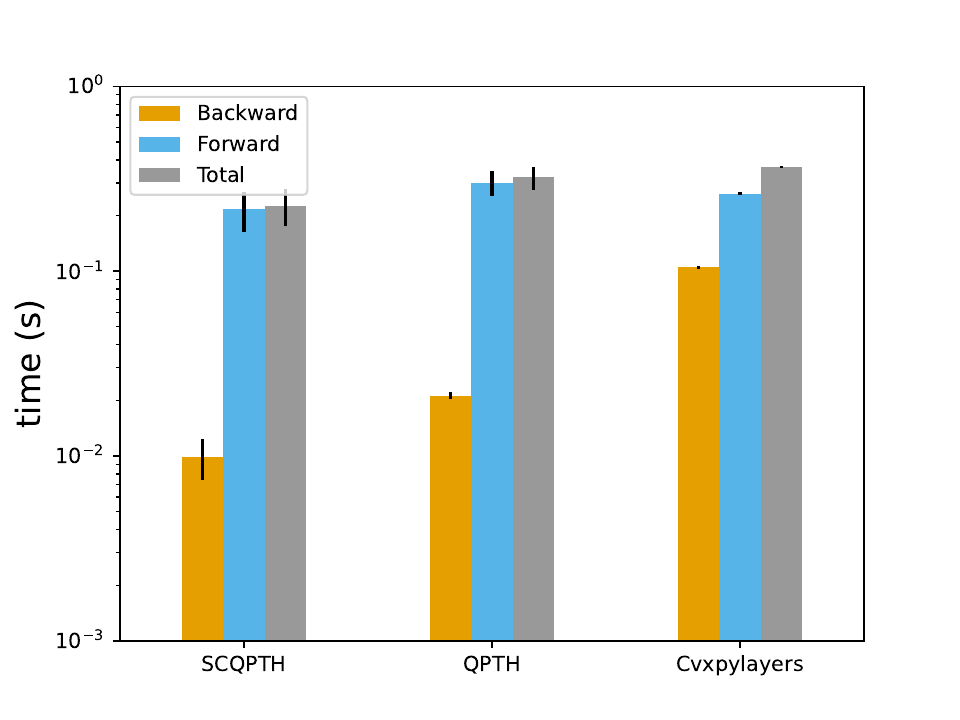}
    \caption{$n= 50, m=100$.}
  \end{subfigure}
  \begin{subfigure}[b]{0.22\linewidth}
    \includegraphics[width=\linewidth , trim={0mm 0cm 0cm 0cm},clip]{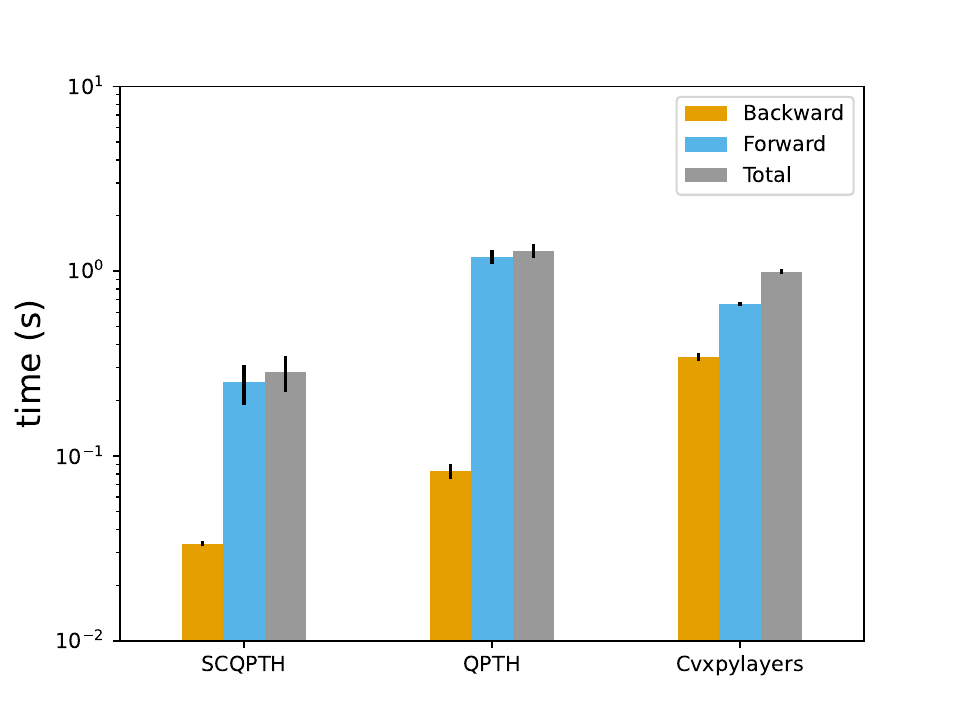}
    \caption{$n = 100, m=200$.}
  \end{subfigure}
  \begin{subfigure}[b]{0.22\linewidth}
    \includegraphics[width=\linewidth , trim={0mm 0cm 0cm 0cm},clip]{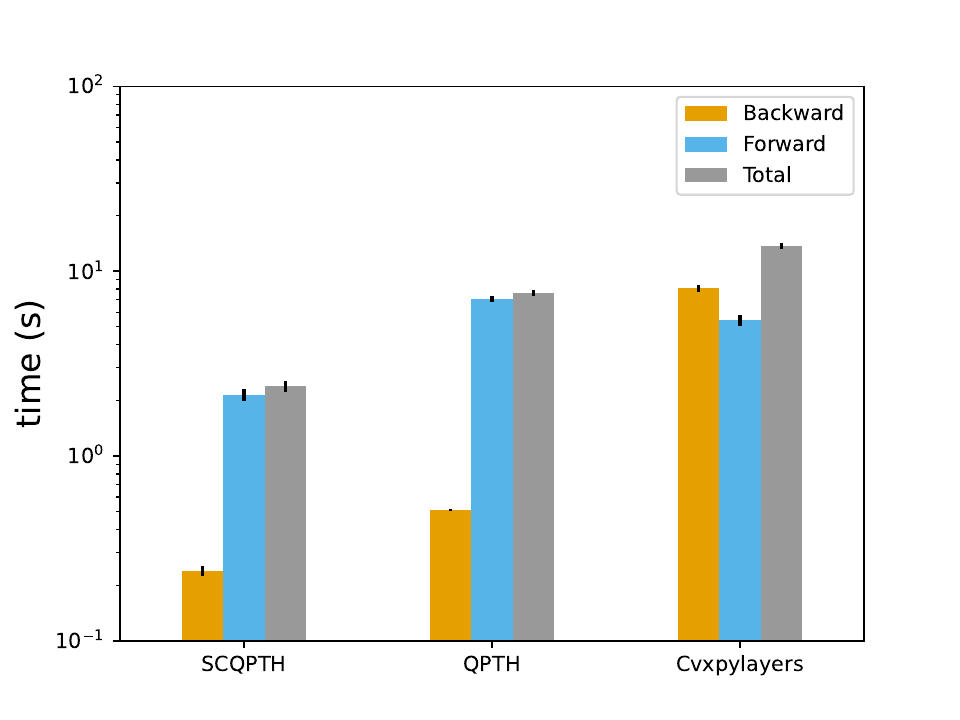}
    \caption{$n = 250, m=500$.}
  \end{subfigure}
  \begin{subfigure}[b]{0.22\linewidth}
    \includegraphics[width=\linewidth , trim={0mm 0cm 0cm 0cm},clip]{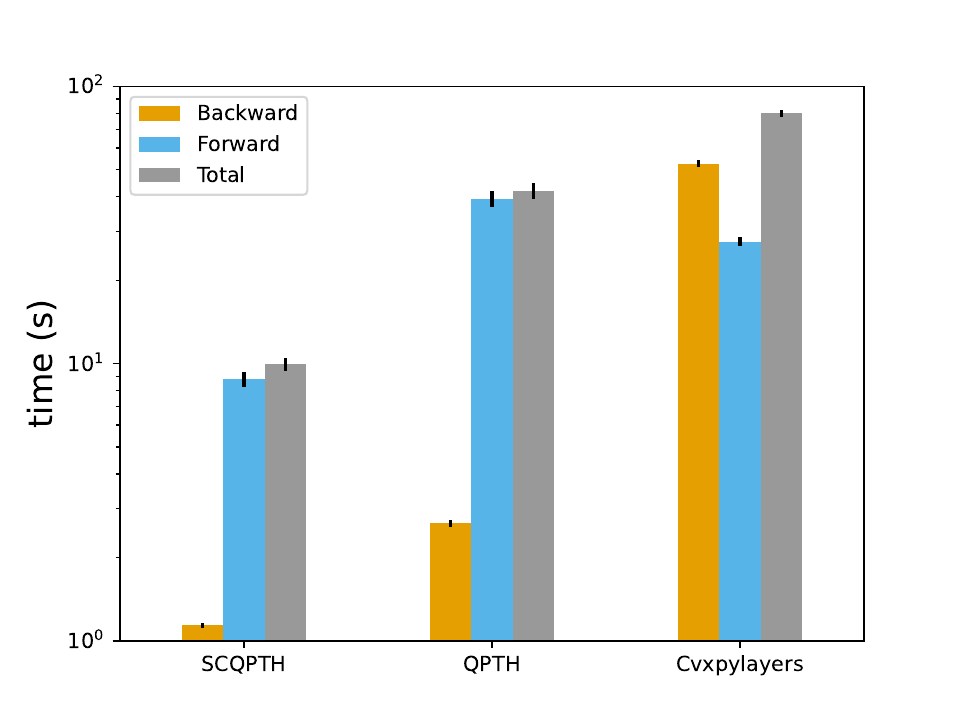}
    \caption{$n = 500, m=1000$.}
  \end{subfigure}
  \begin{subfigure}[b]{0.22\linewidth}
    \includegraphics[width=\linewidth , trim={0mm 0cm 0cm 0cm},clip]{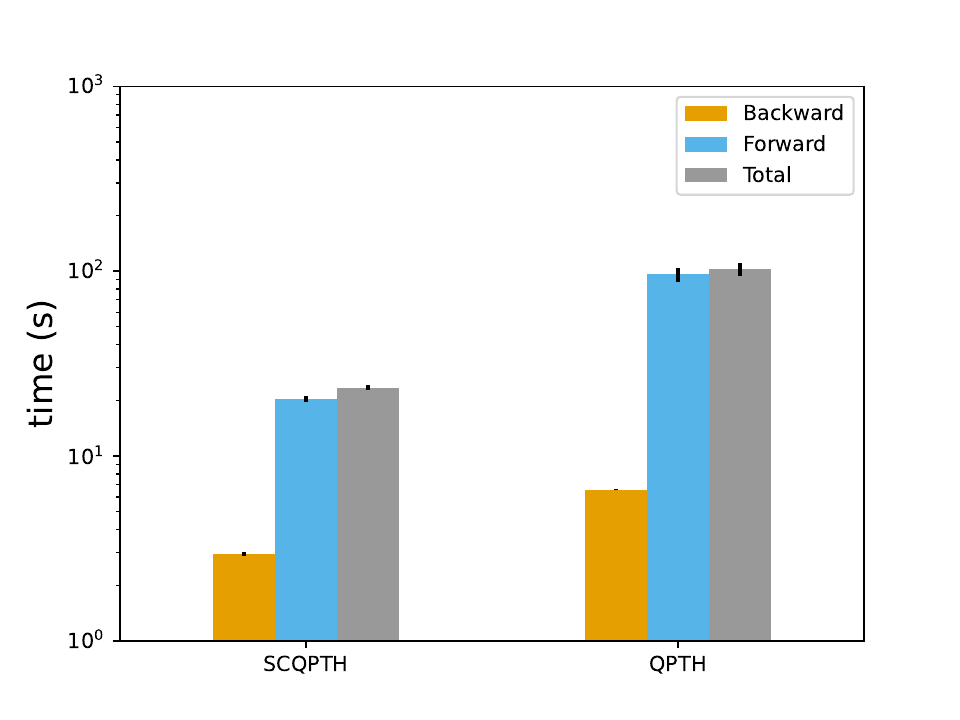}
    \caption{$n = 750, m=1500$.}
  \end{subfigure}
  \begin{subfigure}[b]{0.22\linewidth}
    \includegraphics[width=\linewidth , trim={0mm 0cm 0cm 0cm},clip]{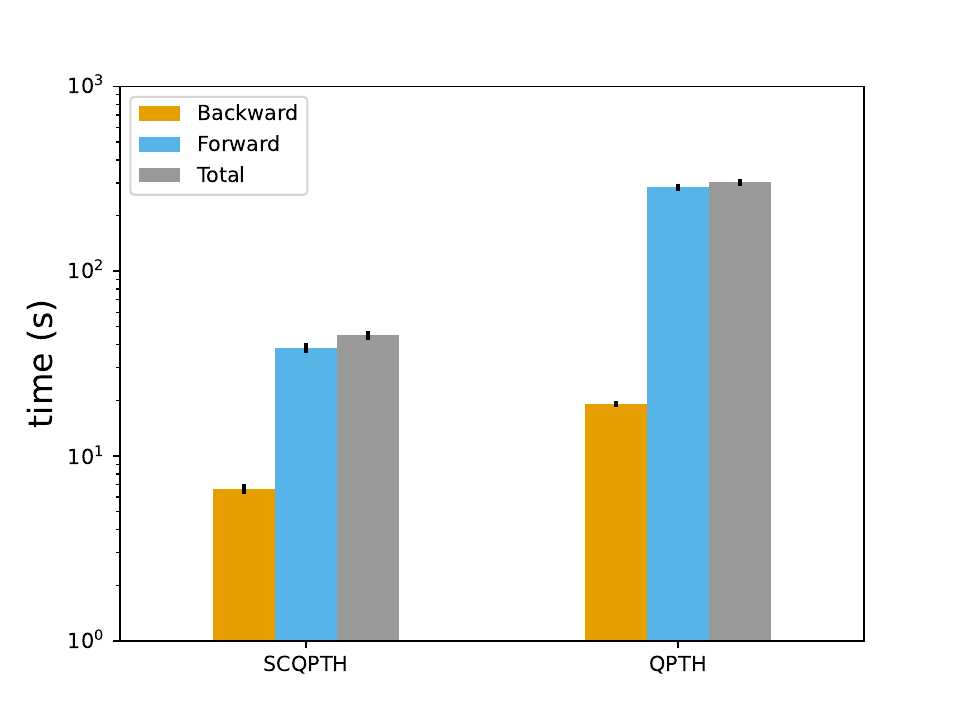}
    \caption{$n = 1000, m=2000$.}
  \end{subfigure}
  \caption{Computational performance of SCQPTH,  QPTH and Cvxpylayers for random QPs of various problem sizes, $n$, constraints $m=2n$, and high stopping tolerance $(1\mathrm{e}{-5})$.}
  \label{fig:exp_2_high_2}
\end{figure}

Finally, Figures \ref{fig:exp_2_low_2} and  \ref{fig:exp_2_high_2}  provide the median runtime and $95\%$-ile confidence interval with a low and high stopping tolerance, respectively, and constraints $m=2n$. With the exception of $n=10$, the SCQPTH layer is the most computationally efficient method and can provide anywhere from a $1\times$ to over $10\times$ improvement in computational efficiency compared to QPTH and Cvxpylayers. For example, when $n=1000$ and stopping tolerance is low, the SCQPTH layer has a median total runtime of $25$ seconds whereas QPTH has median total runtime of $315$ seconds; a $12.6\times$ improvement in efficiency.

\subsection{Experiment 3:  Learning $\bp$}\label{sec:results_3}
We now consider the task of learning a parameterized model for the variable $\bp$:
\begin{equation}
\bp(\btheta) = \bw^T \btheta.
\end{equation}
We follow the procedure outlined in Section \ref{sec:results_2} to generate randomly constrained QPs. We generate random feature variables $\bw \in \mathbb{R}^{5}$ and generate the ground truth coefficients $\btheta^* \in \mathbb{R}^{5 \times n}$ with entries $\btheta^*_{ij} \sim \mathcal{N}(0, 1)$. We set $n=100$, an absolute and relative stopping tolerance of $1\mathrm{e}{-3}$ and consider the case where $m = n$, $m = 2n$ and $m=5n$.  Experiment results are averaged over 10 independent trials. Each trial consists of $100$ epochs,  with total training sample size of 128 and a mini-batch size of 32.

\begin{figure}[H]
  \centering
  \begin{subfigure}[b]{0.35\linewidth}
    \includegraphics[width=\linewidth, trim={0cm 0cm 0cm 0cm},clip]{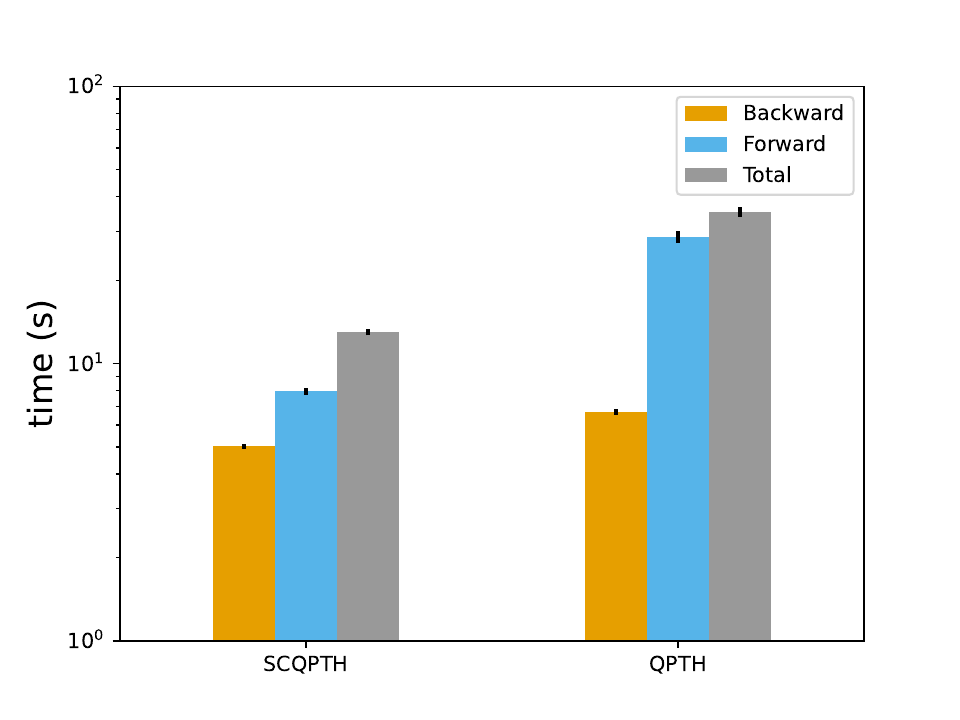}
    \caption{Training Loss.}
  \end{subfigure}
  \begin{subfigure}[b]{0.35\linewidth}
    \includegraphics[width=\linewidth, trim={0cm 0cm 0cm 0cm},clip]{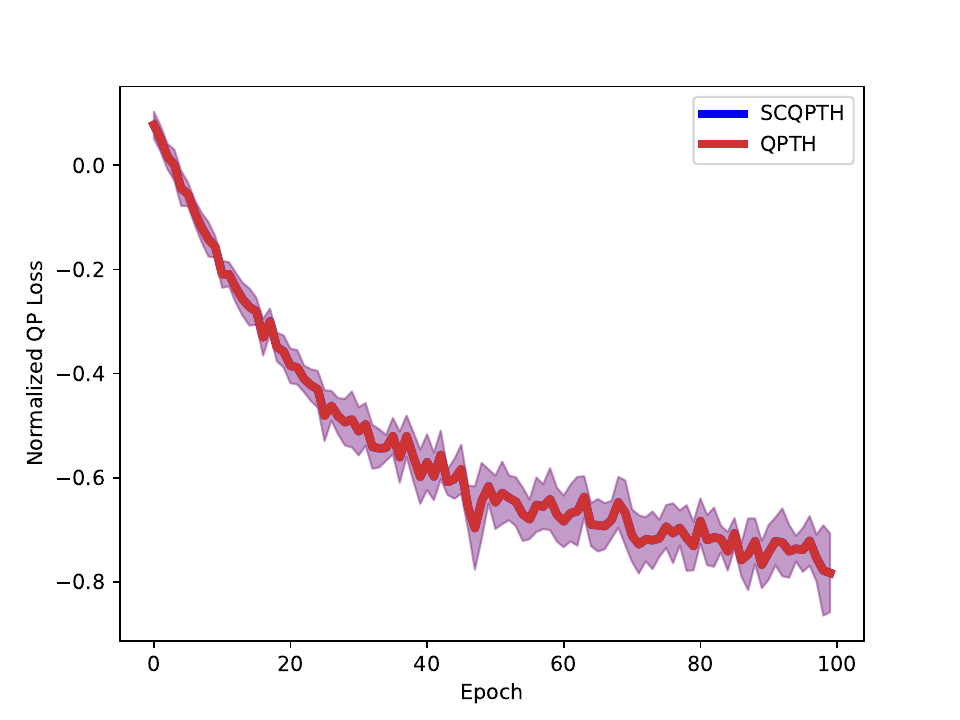}
    \caption{Computational Performance.}
  \end{subfigure}
  \caption{Training loss and computational performance for learning $\bp$ $(n = 100, m = 100)$.}
  \label{fig:learn_p_100}
\end{figure}

\begin{figure}[h]
  \centering
  \begin{subfigure}[b]{0.35\linewidth}
    \includegraphics[width=\linewidth, trim={0cm 0cm 0cm 0cm},clip]{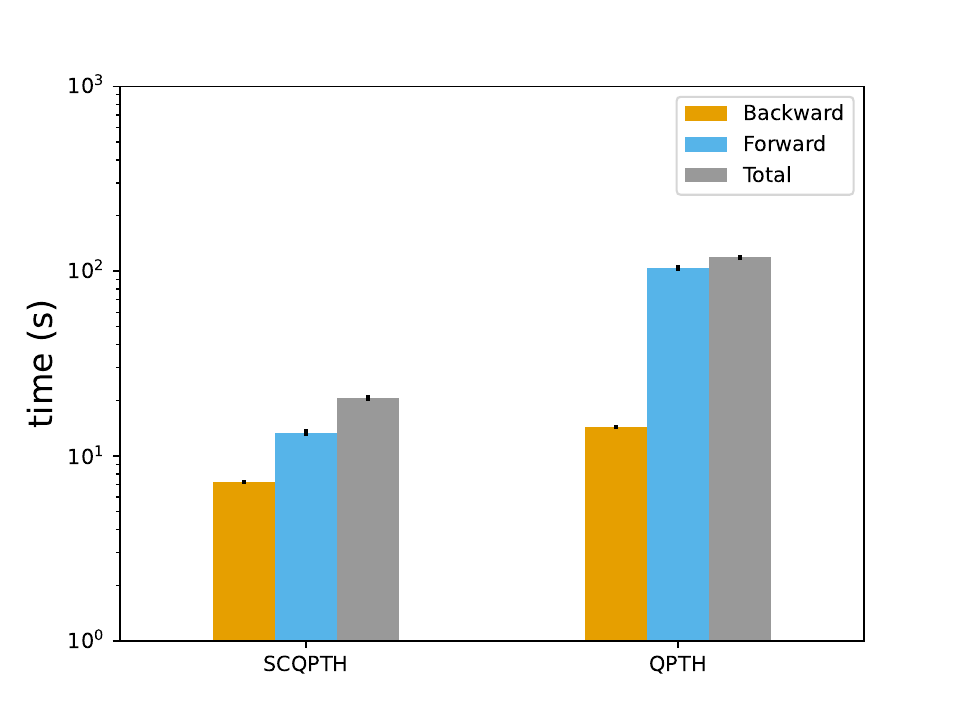}
    \caption{Training Loss.}
  \end{subfigure}
  \begin{subfigure}[b]{0.35\linewidth}
    \includegraphics[width=\linewidth, trim={0cm 0cm 0cm 0cm},clip]{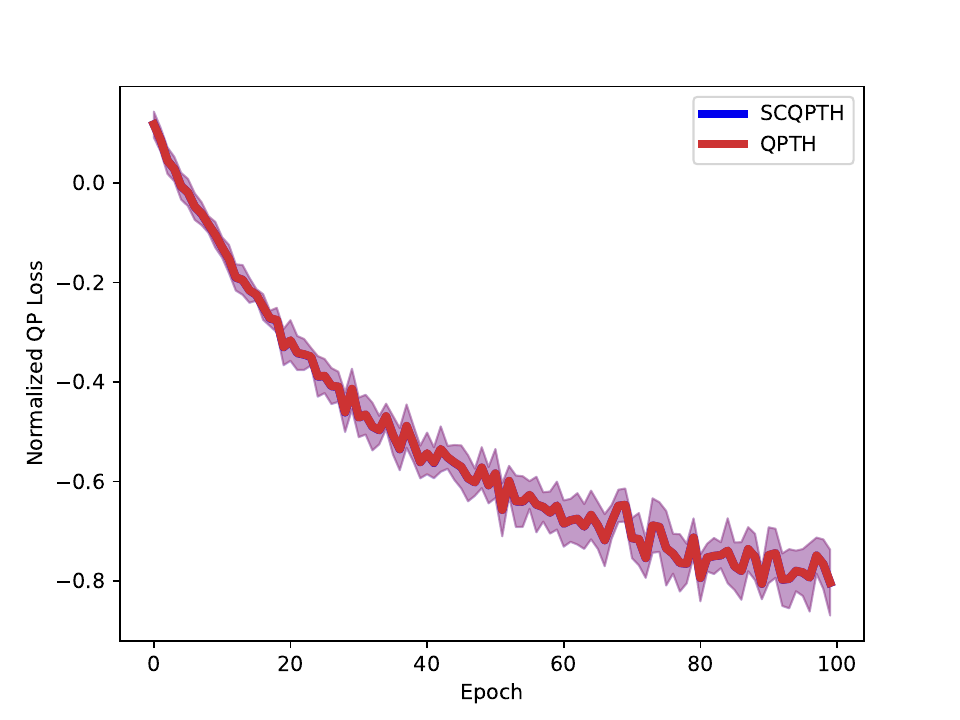}
    \caption{Computational Performance.}
  \end{subfigure}
  \caption{Training loss and computational performance for learning $\bp$ $(n = 100, m = 200)$.}
  \label{fig:learn_p_200}
\end{figure}

\begin{figure}[h]
  \centering
  \begin{subfigure}[b]{0.35\linewidth}
    \includegraphics[width=\linewidth, trim={0cm 0cm 0cm 0cm},clip]{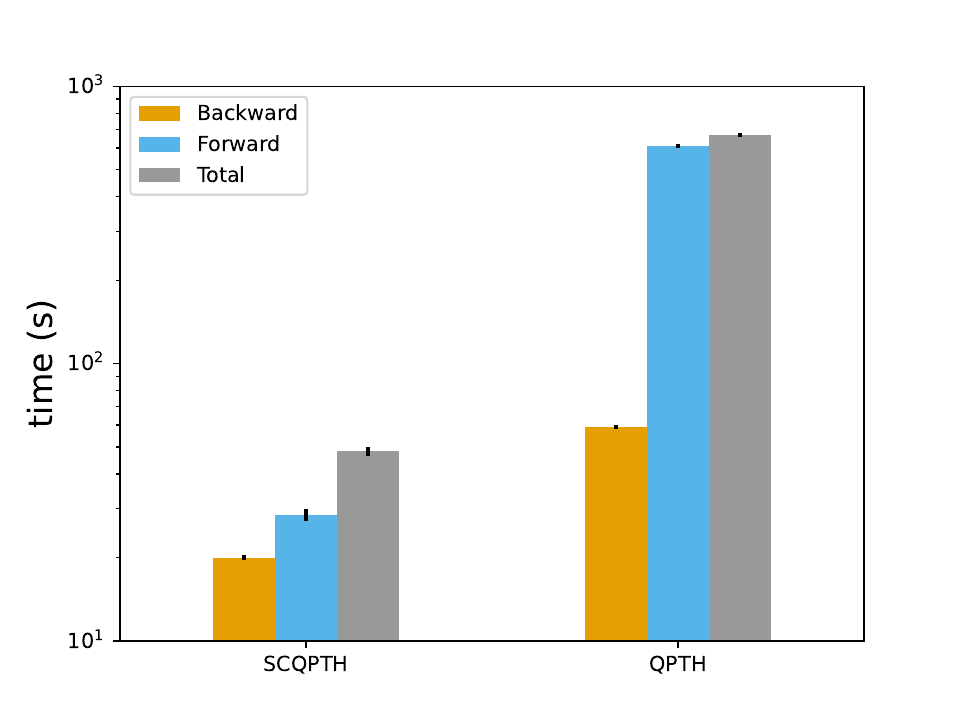}
    \caption{Training Loss.}
  \end{subfigure}
  \begin{subfigure}[b]{0.35\linewidth}
    \includegraphics[width=\linewidth, trim={0cm 0cm 0cm 0cm},clip]{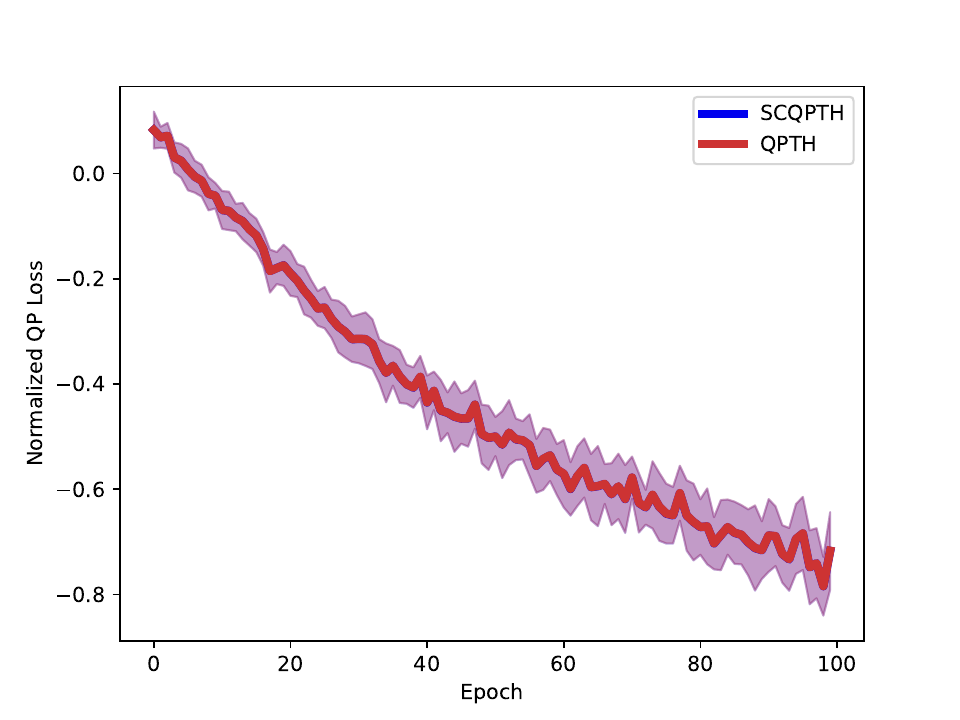}
    \caption{Computational Performance.}
  \end{subfigure}
  \caption{Training loss and computational performance for learning $\bp$ $(n = 100, m = 500)$.}
  \label{fig:learn_p_500}
\end{figure}

We compare the computational efficiency and performance accuracy of SCQPTH with the interior-point solver QPTH.  Figures \ref{fig:learn_p_100}(a) - \ref{fig:learn_p_500}(a) report the average and $95\%$-ile confidence interval training loss at each epoch. We note that the training loss curves are \textit{identical} in all cases; suggesting equivalent training accuracy. Figures \ref{fig:learn_p_100}(b) - \ref{fig:learn_p_500}(b), compares the median runtime to perform $100$ training epochs. When $n = m =100$, SCQPTH  and QPTH requires 13.0 seconds and 35.30 seconds to train; an approximate $2.7\times$ increase in computational efficiency. However, when $n=100$ and $m=500$, the learning process takes approximately $660$ seconds to train the QPTH model, but less than $50$ seconds to train SCQPTH; an over $13\times$ improvement in computational efficiency. 

\section{Conclusion and future work}\label{sec:conclusion}
In this paper, we presented SCQPTH:  a differentiable first-order \textit{splitting} method for \textit{convex quadratic programs}. SCQPTH solves convex QPs using the ADMM algorithm and computes gradients by implicit differentiation of the corresponding fixed-point mapping. Computational experiments demonstrate that for large scale QPs with $100 - 1000$ decision variables, SCQPTH can provide up to an order of magnitude improvement in computational efficiency in comparison to existing differentiable QP methods. Furthermore, in contrast to existing methods,  SCQPTH scales well with the number of constraints and can efficiently handle QPs with thousands of linear constraints.

Our results should be interpreted as a proof-of-concept and we seek to perform further testing on real-world data sets. Moreover, the SCQPTH implementation has several limitations and areas for further improvement. For example, in Section \ref{sec:results_2}, when the number of constraints $m = n$ and the stopping tolerance was high, we observed negligible improvement in computational performance in comparison to the interior-point method. Indeed, it is not uncommon for ADMM and first-order methods in general to have comparatively slow convergence on high accuracy solutions. Accelerating first-order methods \citep{Anderson1965, Sopa2019,Walker2011} has shown to improve the convergence rates in the event that a high accuracy solution is required and is an interesting area of future research. Furthermore, the forward and backward methods of SCQPTH implement several heuristic techniques for preconditioning, scaling, parameter selection and gradient approximation. Providing stronger theoretical justification for the proposed implementations is another area of future research.

\bibliographystyle{plainnat}
\bibliography{Bibliography/Bibliography}

\appendix
\section{Appendix}
\subsection{SCQP Software Implementation}
The SCQPTH software is made available as an open-source Python package, available here:  $$\text{https://github.com/ipo-lab/scqpth}.$$ All benchmark experiments are available here: $ \text{https://github.com/ipo-lab/scqpth\_bench}$. 

The SCQPTH solver interface is implemented as a Python object class which, upon instantiation, takes as input the problem variables: $\bQ$, $\bp$, $\bA$, $\bl$, and $\bu$, and a dictionary of optimization control parameters. The default parameters and a description are presented in Table \ref{table:scqp_params} below.
\begin{table}[h]
\begin{center}
\scriptsize
\begin{tabular}{| l || l ||  l |}
\hline
\textbf{Parameter Name} & \textbf{Description} & \textbf{Default Value} \\ 
 \hline
\textbf{max\_iters} & Maximum number of ADMM iterations. & $1\mathrm{e}{4}$ \\  
\textbf{eps\_abs} & Absolute stopping tolerance. & $1\mathrm{e}{-3}$ \\  
\textbf{eps\_rel} & Relative stopping tolerance. & $1\mathrm{e}{-3}$ \\  
\textbf{eps\_infeas} & Infeasibility stopping tolerance. & $1\mathrm{e}{-4}$  \\  
\textbf{check\_solved} & Interval for checking termination condtions  . & $25$ \\  
\textbf{check\_feasible} & Interval for checking infeasibility conditions. & $25$ \\  
\textbf{alpha} & ADMM relaxation parameter. & $1.2$ \\  
\textbf{alpha\_iter} & Number of initial non-relaxed ADMM iterations. & $100$ \\  
\textbf{rho} & ADMM step-size parameter. & None (auto) \\  
\textbf{rho\_min} & Lower bound on rho. & $1\mathrm{e}{-6}$ \\  
\textbf{rho\_max} & Upper bound on rho. &$1\mathrm{e}{6}$ \\  
\textbf{adaptive\_rho} & True/False for adaptive rho selection. & True \\  
\textbf{adaptive\_rho\_tol} & Threshold for changing rho. & $10$ \\  
\textbf{adaptive\_rho\_iter} & Number of initial non-adaptive ADMM iterations.. & $50$ \\  
\textbf{adaptive\_rho\_max\_iter} & Maximum number of adaptive iterations. & $1\mathrm{e}{3}$ \\ 
\textbf{sigma} & Tikhonov regularization for semi-definite problems. & $0.0$ \\   
\textbf{scale} & True/False for automatic scale. & True \\  
\textbf{beta} & Shrinkage factor for automatic scaling. & None (auto) \\  

 \hline
\end{tabular}
\end{center}
\caption{\label{table:scqp_params} SCQPTH software optimization control parameters and default values.}
\end{table}
\subsection{Proof of Proposition $\ref{prop:admm_fixed_point}$}\label{app:prop_1}
We define $\bv^k = \bA \bx^{k+1} +\bmu^k$. We can therefore express Equation $\eqref{eq:admm_qp_z}$ as:
\begin{equation}\label{eq:app_z}
\bz^{k+1} =  \Pi( \bA \bx^{k+1}  + {\bmu}^k )  = \Pi( \bv^k ),
\end{equation}
and Equation $\eqref{eq:admm_qp_mu}$ as:
\begin{equation}\label{eq:app_mu}
{\bmu}^{k+1} = {\bmu}^{k} +  \bA \bx^{k+1} - \bz^{k+1} = \bv^k - \Pi( \bv^k ).
\end{equation}
Substituting Equations $\eqref{eq:app_z}$ and $\eqref{eq:app_mu}$ into Equation $\eqref{eq:admm_qp_x}$ gives the desired fixed-point iteration:

\begin{equation}\label{eq:app_admm_fp}
\bv^{k+1} =  \bA \bx^{k+1} +\bmu^k = \bA [\bQ + \rho \bA^T\bA]^{-1} (- \bp + \rho \bA^T (2\Pi(\bv^{k}) - \bv^{k})) + \bv^k - \Pi(\bv^{k})
\end{equation}

\subsection{Proof of Proposition $\ref{prop:admm_grads}$}\label{app:prop_2}

We define $F \colon \mathbb{R}^m \to \mathbb{R}^m$ as:
\begin{equation}\label{eq:app_F}
F(\bv) = \bA [\bQ + \rho \bA^T\bA]^{-1} (- \bp + \rho \bA^T (2\Pi(\bv) - \bv)) + \bv - \Pi(\bv)
\end{equation}
and let
\begin{equation}\label{eq:app_M}
\bM = [\bQ + \rho \bA^T\bA] \qquad \text{and} \qquad \blr = - \bp + \rho \bA^T (2\Pi(\bv) - \bv)
\end{equation}
Taking the partial differentials of Equation $\eqref{eq:app_F}$ with respect to the relevant problem variables therefore gives:
\begin{equation}\label{eq:app_partial_F}
\begin{split}
\partial F(\bv,\betta)  & = \partial \bA \bM^{-1} \blr - \bA \bM^{-1} \partial \bM \bM^{-1} \blr +  \bA \bM^{-1} \partial \blr\\
& = \partial \bA \bx^* - \bA \bM^{-1} \partial \bM \bx^* + \bA \bM^{-1} \partial \blr\\
& = \partial \bA \bx^* - \bA \bM^{-1}\Big( \frac{1}{2} (\partial \bQ + \partial \bQ^T) + \partial \bp + \rho \partial A^T \bmu^* \Big)
\end{split}
\end{equation}

Substituting the gradient action of Equation $\eqref{eq:app_partial_F}$ into Equation $\eqref{eq:jacob_F}$ and taking the left matrix-vector product of the transposed Jacobian with the previous backward-pass gradient, $\frac{\partial \ell }{\partial \bx^*}$, gives the desired result.

\subsection{Proof of Proposition $\ref{prop:admm_grads_box}$}\label{app:prop_3}
Taking the partial derivatives of the KKT optimality conditions $\eqref{eq:qp_kkt}$ gives the following linear system of equations:
\begin{equation} \label{eq:diff}
\begin{split}
\begin{bmatrix}
 \bQ &  \bA^T &   \bA^T \\
\diag(\by_-^*)\bA & \diag (\bA \bx^* - \bl ) & 0\\
\diag(\by_+^*)\bA & 0 & \diag (\bA \bx^* - \bu ) 
\end{bmatrix}
\begin{bmatrix}
{\partial \bx }\\
{\partial \by_- }\\
{\partial \by_+}
\end{bmatrix}
= -
\begin{bmatrix}
{\partial \bQ} \bx^* + { \partial \bp } + {\partial  \bA}^T \by_-^* + {\partial  \bA}^T \by_+^* \\
\diag(\by_-^*) {\partial  \bA } \bx^* - \diag(\by_-) \partial  \bl\\
\diag(\by_+^*) {\partial  \bA } \bx^* - \diag(\by_+) \partial  \bu
\end{bmatrix}.
\end{split}
\end{equation}
Implicit differentiation of $\eqref{eq:diff}$ gives the following linear system of equations:
\begin{equation} \label{eq:diff_sol}
\begin{split}
\begin{bmatrix}
\bd_{\bx} \\
\bd_{\by_-} \\
\bd_{\by_+}
\end{bmatrix}
= -
\begin{bmatrix}
\bQ &  \bA^T \diag(\by_-^*)&   \bA^T \diag(\by_+^*) \\
\bA & \diag (\bA \bx^* - \bl ) & 0\\
\bA & 0 & \diag (\bA \bx^* - \bu ) 
\end{bmatrix} ^{-1}
\begin{bmatrix}
\big( \frac{\partial \ell }{\partial \bz^*} \big)^T \\
0\\
0
\end{bmatrix}.
\end{split}
\end{equation}

Let $\bd_{\by} = \diag(\by_-) \bd_{\by_-} + \diag(\by_+) \bd_{\by_+}$. Then from equation $\eqref{eq:diff_sol}$ we have:
\begin{equation}\label{eq:app_d_lambda}
\bA^T\bd_{\by} = \Big( - \Big(\frac{\partial \ell }{\partial \bz^*}\Big)^T - \bQ \bd _{\bx} \Big).
\end{equation}
The pseudo-inverse of $\bA^T$ therefore approximates a solution $\bd_{\by}$; given by:
\begin{equation}
\bd_{\by} =  [\bA^T]^\dagger \Big( - \Big(\frac{\partial \ell }{\partial \bx^*} \Big)^T - \bQ  \bd_{\bx} \Big),
\end{equation}

 Combining the complimentary slackness condition with Equation $\eqref{eq:diff_sol}$ therefore uniquely determines the relevant non-zero elements of $\bd_{\by_-}$ and $\bd_{\by_+}$. Specifically:
\begin{equation}\label{eq:dy_neg}
\bd_{\by_{-_j}}  = \begin{cases}
                \bd_{\by_j} / \by_{-_j} & \text{if }   \by_{-_j} \leq 0\\
                0 & \text{otherwise,} \\
                \end{cases} \qquad \text{and} \qquad
\bd_{\by_{+_j}} =  \begin{cases}
                \bd_{\by_j} / \by_{+_j} & \text{if }   \by_{+_j} \geq 0\\
                0 & \text{otherwise,} \\
                \end{cases}
\end{equation}
 Computing the  left  matrix-vector  product  of  the  Jacobian  with  the  previous  backward-pass  gradient, $\frac{\partial \ell }{\partial \bx^*}$ gives the desired gradients.

\end{document}